\documentclass[12pt]{amsart}
\usepackage{amsmath}
\usepackage{MnSymbol}
\usepackage{amsfonts}\usepackage{verbatim}\usepackage{bbold}
\usepackage{ifsym}\usepackage{marvosym}\usepackage{fontawesome}
\usepackage{dictsym}\usepackage{wasysym}
\usepackage{amstext}
\usepackage{amsbsy}
\usepackage{amsopn}\usepackage{hyperref}
\usepackage{amsthm}\usepackage{color}
\DeclareMathSymbol{\varkappa}       {\mathord}{AMSb}{"7B}

\def\proclaim#1{\vskip0.5em\noindent{\bf #1}\it }
\def\endproclaim{\vskip0.5em\par\noindent\rm}

\def\proclaim#1{\vskip0.5em\noindent{\bf #1}\it}
\def\endproclaim{\vskip0.5em\par\noindent\rm}

\def\undersetbrace#1\to#2{\underbrace{#2}_{#1}}

\def\demo#1{\vskip0.5em\noindent{\bf #1\ }}

\def\text#1{\mbox{#1}}
\def\flushpar{\par\noindent}

\def\mod{\mbox{ mod }}

\newcommand{\mapright}[1]{%
    \smash{\mathop{%
        \hbox to 1cm{\rightarrowfill}
        }
    \limits^{#1}
    }
}
\newcommand{\mapleft}[1]{%
    \smash{\mathop{%
        \hbox to 1cm{\rightarrowfill}
        }
    \limits_{#1}
    }
}

\def\e{\epsilon}
\def\a{\alpha}
\def\b{\beta}
\def\G{\Gamma}
\def\g{\gamma}
\def\d{\delta}
\def\D{\Delta}
\def\s{\sigma}

\def\th{\theta}
\def\L{\Lambda}\def\l{\lambda}

\def\x{\times}

\def \C{\mathcal C}
\def\f{\flushpar}
\def\u{\underline}
\def\v{\varphi}

\def\om{\omega}
\def\Om{\Omega}
\def\B{\mathcal B}
\def\T{\widehat T}
\def\({\biggl(}
\def\){\biggr)}
\def\<{\langle}
\def\>{\rangle}
\def\bdy{\partial}
\def\bul{\smallskip\f$\bullet\ \ \ $}\def\lfl{\lfloor}\def\rfl{\rfloor}

\def\){\biggr)}

\def\rv{\text{random variable}}

\def\with respect to{\text{with respect to}}

\def\<{\bold\langle}
\def\>{\bold\rangle}

\def\bul{\smallskip\f$\bullet\ \ \ $}\def\sms{\smallskip\f}\def\sbul{\f$\bullet\
\ \ $}\def\Lra{\Longrightarrow}\def\pf{\smallskip\f{\it Proof} }\def\Par{\smallskip\f\P}
\def\Smi{\smallskip\f{\LARGE \smiley}}
\def\ycnu{(Y,\nu)}

\begin{document}

\title{ Local limit theorems for suspended semiflows }
\author{Jon. Aaronson $\&$ Dalia Terhesiu}
 \address[Aaronson]{\ \ School of Math. Sciences, Tel Aviv University, 69978 Tel Aviv, Israel.\f\ \ \  {\it Webpage }: 
 {\tt http://www.math.tau.ac.il/$\sim$aaro}}\email{aaro@tauex.tau.ac.il}
\

\address[Terhesiu]{Department of Mathematics, University of Exeter, Exeter EX4 4QF, UK}\email[Terhesiu]{D.Terhesiu@exeter.ac.uk}
\subjclass[2010]{37A40, 60F05}
\keywords{Local limit  mixing, rational weak mixing, operator local limit theorem, suspended semiflow,  cocycle, continuous time random walk,  
  fibered system, transfer operator, quasicompact action, {\tt G-M} map, AFU map,
 infinite ergodic theory.}
\thanks { Aaronson's research was partially supported by ISF grant No. 1289/17.}

\begin{abstract}  We prove   local limit theorems for a cocycle over a semiflow by establishing topological, 
mixing properties of the 
associated skew product semiflow.
We also establish conditional rational weak mixing of certain 
  skew product semiflows and various mixing properties including order 2 rational weak mixing of hyperbolic geodesic flows on the tangent spaces of cyclic covers. 19/6/2020
 \end{abstract}\dedicatory{Dedicated to the memory of Nat Friedman}
\maketitle\markboth{Local limits via mixing}{Jon. Aaronson $\&$ Dalia Terhesiu}
\section*{\S0 Introduction}
\subsection*{Semiflows}
\

We denote a standard $\s$-finite, measure space by $(X,m)$. The space $X$ is Polish (a complete, separable metric space) and the measure $m$ is defined on $\B=\B(X)$ the Borel subsets of $X$.   
\

Let
$$\text{\tt MPT}\,(X,m):=\{\text{measure preserving transformations of}\, (X,m)\}.$$

The collection $\text{\tt MPT}\,(X,m)$ is a   Polish semigroup under composition with respect to the {\tt weak operator topology} defined by $T_n\to T$ if
$$1_A\circ T_n\xrightarrow[n\to\infty]{m}1_A\circ T\ \forall\ A\in\B.$$
\

A {\it semiflow}  on $(X,m)$ is a continuous (semigroup) homomorphism $\Phi:\Bbb R_+\to\text{\tt MPT}\,(X,m)$. 
\

A  {\it flow}  on $(X,m)$ is a continuous (group) homomorphism $\Phi:\Bbb R\to\text{\tt MPT}\,(X,m)$. 

\

The semiflow  $\Phi:\Bbb R_+\to\text{\tt MPT}\,(X,m)$ is called {\it invertible} if  each transformation $\Phi_t$ is invertible.
In this case,  $\Phi$ embeds in a flow 
$\Psi:\Bbb R\to\text{\tt MPT}\,(X,m)$ defined by
$\Psi_t:=\Phi_{|t|}^{\text{\tiny sgn}(t)}\ \ (t\ne 0)\ \&\ \Psi_0:=\text{\tt Id}$.

\

\subsection*{Cocycles}
\

Let $\Bbb G$ be a locally compact, Polish Abelian group and let $(X,m,T)$ be a  measure preserving transformation.
\

Given a  measurable function $c:X\to\Bbb G$, we may define 
$C:\Bbb N\x X\to\Bbb G$ by  $C(n,x)=c_n(x):=\sum_{k=0}^{n-1}c(T^kx)$ and obtain the {\tt cocycle equation} $C(n+n',x)=C(n,x)+C(n',T^n(x))$. 
\

Analogously, for $(X,m,\Phi)$  a  semiflow,  a  $\Bbb G$-valued  $\Phi$-{\it cocycle} is a  measurable function 
$C:\Bbb R_+\x X\to\Bbb G$ satisfying $C(t+t',x)=C(t,x)+C(t',\Phi_t(x))$.
\

A $\Phi$-cocycle is a $\Phi$-{\it coboundary} if it has the form $(t,x)\mapsto H(x)-H(\Phi_t(x))$ where $H:X\to\Bbb G$ is measurable and two $\Phi$-cocycles
are {\it cohomologous} if they differ by a $\Phi$-coboundary.

\subsection*{Local limit theorems}
In this paper we study the local limit theory of cocycles over a semiflow.
Let $\Bbb G\le\Bbb R^\kappa$ be a closed subgroup of full dimension, let  $(X,m,\Phi)$ be a probability preserving (abbr. {\tt PP})
semiflow and let
$C:\Bbb R_+\x X\to\Bbb G$ be a $\Phi$-{cocycle}. For $0<p\le 2$, let $\mathcal S$ be a globally supported, strictly $p$-stable  random variable
on $\Bbb R^\kappa$ and let $b:\Bbb R_+\to\Bbb R_+$ be $\tfrac1p$-regularly varying. Fix a ring $\mathcal R\subset\B(X)$.
\

The $\Phi$-{cocycle} $C:\Bbb R_+\x X\to\Bbb G$ is said to satisfy:
\

-- the {\it integrated $(\mathcal S,b)$-local limit theorem} (with respect to $\mathcal R$)
if
\begin{align*}\tag{{\tt Int-LLT}}\label{Int-LLT} 
b(t)^{\kappa}m(A&\cap \Phi_t^{-1}B\cap [C(t,\cdot)\in x(t)+U])\\ &\xrightarrow[t\to\infty,\ \frac{x(t)}{b(t)}\to z]{} 
f_\mathcal S(z)\cdot m(A)m(B)m_\Bbb G(U)\\ & \forall\ A,B\in\mathcal R\ \&\ U\subset\Bbb G\  \text{precompact},\ m_\Bbb G(\bdy U)=0;
\end{align*} where $f_\mathcal S:\Bbb R^d\to\ [0,\infty)$ is the probability density function of $\mathcal S$ and 
\

-- the {\it conditional $(\mathcal S,b)$-local limit theorem} (with respect to $\mathcal R$) if
\begin{align*}\tag{{\tt Con-LLT}}\label{Con-LLT}  
b(t)^{\kappa}\widehat{\Phi_{t}}(&1_{A\cap [C(t,\cdot)\in x(t)+U]})\xrightarrow[t\to\infty,\ \frac{x(t)}{b(t)}\to z]{} 
f_\mathcal S(z)\cdot m(A)m_\Bbb G(U)\\ & \text{a.e.}\ \forall\ A\in\mathcal R\ \&\ U\subset\Bbb G\  \text{precompact},\ m_\Bbb G(\bdy U)=0.\end{align*} 
Here and throughout, for $T\in\text{\tt MPT}\,(X,m)$,  ${\T}:L^1(m)\to L^1(m)$ denotes the 
  {\it transfer operator} of $T$   defined by
$$\int_X{\T}f\cdot gdm=\int_Xf\cdot g\circ Tdm.$$

The convergence (in both ({\tt Con-LLT}) and ({\tt Int-LLT})) is often 
uniform in $z$ in compact subsets of $\Bbb G$, and the convergence in  ({\tt Con-LLT}) is often also uniform on subsets in the ring $\mathcal R$.
\

Versions of ({\tt Int-LLT}) in the case where $p=2$ (i.e. $\mathcal S$ is Gaussian) can be found in
\cite{Waddington},\ \cite{Iwata} and \cite{DN}.

\subsection*{What's new}\ \ 
\

In this paper,  we give conditions for ({\tt Con-LLT})   and for ({\tt Int-LLT})
 (theorem 1). 
\

Our approach is via mixing properties of the associated skew product semiflow. 
 \
 
 Conditions are obtained for {\tt conditional rational weak mixing} of skew product semiflows (theorem 2).
 \
 
 Applying all this, we  establish  various mixing properties   for the geodesic flow of a cyclic cover of a compact hyperbolic surface
 (theorem 3).
  
\

\section*{\S1 Background}

\subsection*{ Skew product semiflows}
    
Let  $(Y,\nu,\Phi)$ be  a measure preserving semiflow and let $G:\Bbb R_+\x Y\to\Bbb G\le\Bbb R^\kappa$
(a closed subgroup of full dimension) be  a $\Phi$-cocycle.
\

The {\it skew product semiflow} over the {\it base}  $(Y,\nu,\Phi)$  with {\it displacement} $G$ is the measure preserving semiflow 
 $(Y\x\Bbb G,\nu\x m_{\Bbb G},\Phi^{(G)})$ defined by
$$\Phi^{(G)}_t(y,z):=(\Phi_t(y),z+G(t,y))\ \ (t>0).$$
The {\it deck holonomies} are defined by
\begin{align*}\tag*{}
 \label{deck}Q_h(y,z):=(y,z+h)\ \ \ (h\in\Bbb G)
\end{align*}
and commute with the semiflow: $Q_h\circ\Phi^{(G)}_t=\Phi^{(G)}_t\circ Q_h$.
\

\subsection*{Suspended semiflow and renewal process }
\

Let $(X,m,T)$ be an ergodic {\tt MPT} ({\tt EMPT}).
Let $\frak r:X\to\Bbb R_+$ be measurable (aka the {\it roof function}). 
\

The {\it renewal process} generated is 
$$N(y)(x):=\#\,(\G(x)\cap (0,y])\ \text{where}\ \G(x):=\{\frak r_n(x):=\sum_{k=0}^{n-1}\frak r(T^kx):\ n\ge 1\}.$$
\

In case $m(X)=1$ and $E(\frak r)<\infty$, 
$\tfrac{E(N(t))}t\xrightarrow[t\to\infty]{}\tfrac1{E(\frak r)}$ by the ergodic theorem. The limit
$\l:=\tfrac1{E(\frak r)}>0$ is called
the {\it asymptotic intensity} of the renewal process.
\

The 
{\it suspended semiflow} with {\it section} $(X,m,T)$ is the  semiflow
$(Y,\nu,\Phi)$ where
\begin{align*}
 &Y=\{(x,y)\in X\x\Bbb R_+:\ 0\le y<\frak r(x)\},\\ &    \nu(A\x C)=m(A)\text{\tt Leb}\,(C) \ \ (A\x C\in\B(Y)),
 \\ &
 \Phi_t(x,y)=(T^nx,y+t-\frak r_n(x))\\ &\  \text{where  $n=\mathcal N(t)(x,y)$ is such that}
 \\ & \frak r_n(x)\le y+t<\frak r_{n+1}(x).
\end{align*}
$$\text{\tt i.e.}\ \ \ \ \mathcal N(t)(x,y)=\#\,(\G(x)\cap (0,y+t])=N(t+y)(x).$$
\

It follows from the definitions that $\nu(Y)=\int_X\frak r dm$ and that for $t>0$,
$$t\mapsto\Phi_t\in\text{\tt MPT}\ycnu$$
is continuous.
\

Next it is routine to check that

\begin{align*}\tag{\phone}\label{phone}
 \mathcal N(t+t')=\mathcal N(t)+\mathcal N(t')\circ\Phi_t\ \forall\  (x,y)\in Y,\ t,t'>0;
\end{align*}

whence $\Phi_{t+t'}=\Phi_t\circ\Phi_{t'}$, showing that $\Phi:\Bbb R\to\text{\tt MPT}(Y,\nu)$ 
is indeed a semiflow and $\mathcal N:\Bbb R_+\x Y\to\Bbb R$ is a $\Phi$-{\it cocycle}.
\

\

Accordingly, we'll denote the 
 suspended semiflow with  section $(X,m,T)$  by $(Y,\nu,\Phi)=:(X,m,T)^{\frak r}$ and
 call $\mathcal N:\Bbb R_+\x Y\to\Bbb R$ the {\it  renewal cocycle}. 
\

Note that $Y\cong X\x\Bbb R_+/\sim$ with $(x,y)\sim (x',y')$ iff $\exists
n\ge 0$ so that either  $T_\frak r^n(x,y)=(x',y')$ or $T_\frak r^n(x',y')=(x,y)$ where $T_\frak r(x,y):=(Tx,y-\frak r(x))$.
\subsection*{Topological suspended flows}\label{topflow} \ Let $T$ be a  bi-Lipschitz, homeomorphism of a
Polish space $(X,d)$ and let
$\frak r:X\to\Bbb R_+$ be uniformly continuous and bounded away from $\{0,\infty\}$. 
\

There is a metric  on $Y$ so that $Y$ equipped with the topology generated is   locally homeomorphic 
with neighbourhoods in $X\x\Bbb R_+$, whence $Y$ is Polish; and so that
$\Phi$ is a continuous flow on $Y$. Moreover, if $X$ is [{\tt\small locally}] compact then so is $Y$. See \cite[\S4]{BowenWalters}.

\subsection*{Jump cocycles over a suspended semiflow}
\

Let $(Y,\nu,\Phi)=(X,m,T)^{\frak r}$ be a  suspended semiflow as above, 
let $\Bbb G\le\Bbb R^{\kappa}$ be a subgroup of full dimension
and let $\varphi:X\ \to\Bbb G$ be measurable.
\

Now define the {\it jump cocycle} $J=J^{(\varphi)}:\Bbb R_+\x Y\to\Bbb G$  by
\begin{align*}\tag{\Football}\label{Football}
J(t,(x,y)):=\varphi_{\mathcal N(t)(x,y)}(x)
\end{align*}
where $\varphi_N:=\sum_{k=0}^{N-1}\varphi\circ T^k$.
\

By (\phone), $J$ is a $\Phi$-{\it cocycle}.
\

We refer to both $\varphi:X\ \to\Bbb G$  and  $J^{(\varphi)}:\Bbb R_+\x Y\to\Bbb G$ as the {\it displacement}.

\

By a simple change of variables, the skew product semiflow of the jump cocycle  over a suspended semiflow
is isomorphic with the suspended semiflow over the displacement  skew product.
That is: if $(Y,\nu,\Phi)=(X,m,T)^\frak r$ and  $J:=J^{(\varphi)}$, then
\begin{align*}\tag{\faLinux}\label{faLinux}
 (X\x\Bbb G,\nu\x m_\Bbb G,\Phi^{(J)}) \cong (X\x\Bbb G,m\x m_\Bbb G,T_\varphi)^{\overline{\frak r}}
\end{align*}

where $\overline{\frak r}(x,z):=\frak r(x)$ and $T_\varphi(x,y):=(Tx,y+\varphi(x))$. See \cite{nice}. 
\

 \subsection*{Example: Continuous time random walks}

\

  A {\it  random walk}  ({\tt RW}) on the locally compact, Polish, Abelian group $\Bbb G$ is a sequence of $\Bbb G$-valued 
  $\rv$s $(S_n:\ n\ge 1)$ defined by
  $S_n:=\sum_{k=1}^n X_k$
  with $(X_n:\ n\ge 1)$ independent, identically distributed random variables on $\Bbb G$. Their mutual distribution $f\in\mathcal P(\Bbb G)$ is known as the {\it jump distribution}
  of the {\tt RW}.

\

 A {\it continuous time random walk} (abbr. {\tt CTRW}) (in \cite{Montroll-Weiss},\ {\it jump process} in \cite{Breiman}) on $\Bbb G$   is a one-parameter family
 $(S^{(t)}:\ t>0)$ of $\rv$s on $\Bbb G$ defined by
 $S^{(t)}:=S_{N(t)}$
 where $S_n=\sum_{k=1}^nX_k$ is as above and $N(t):=\#\,\G\cap (0,t]\ \ (t>0)$ where $\G$ is a Poisson point process on $\Bbb R_+$ which is independent of 
 $(X_n:\ n\ge 1)$ and where $E(N(t))=\l t\ \forall\ t>0$ where $\l>0$ (the {\it intensity} of the process).
 \

We'll exhibit a jump cocycle over a suspended semiflow with the same  distribution as $(S^{(t)}:\ t>0)$\label{CTRW-jump}.
\

Let $(X,m,T)$ be given by $X:=(\Bbb G\x\Bbb R_+)^\Bbb N$,   $T=\text{shift}$ and
$m=\mu^\Bbb N$ where
$$\mu(U\x A):=\frac{f(U)}\l\int_Ae^{-\frac{x}\l}dx.$$
Denote $x\in X$ by $x=(\u k,\u r)=((k_1,k_2,\dots),(r_1,r_2,\dots))$ and define 
$\frak r:X\to\Bbb R_+$ by $\frak r(x):=r_1$ and
 $\varphi:X\to\Bbb G$ by $\varphi(x):=k_1$.
\

Let $(X,m,T)^\frak r=(Y,\nu,\Phi)$ and define $J^{(\varphi)}:\Bbb R_+\x Y\to\Bbb G$ as above.
\

We claim that the $\nu$-distribution of $(J^{(\varphi)}(t,\cdot):\ t>0)$ coincides with the distribution of
$(S^{(t)}:\ t>0)$.
\

Indeed, considered with respect to $m,\ \{\frak r_n(\cdot):\ n\ge 1\}=:\G(\cdot)$ is a Poisson point process on $\Bbb R_+$ with 
intensity $\l$. For $y>0$, so is $\G_y(\cdot):=(\G(\cdot)-y)\cap\Bbb R_+$.
\

For   $(x,y)\in Y$, let
$\mathcal N(t)(x,y):=\#(\G(\cdot)\cap (y,y+t])=\#(\G_y(\cdot)\cap (0,t])$, then the distribution
of $(\mathcal N(t)(\cdot,y):\ t>0)$ does not depend on $y>0$. 
\

Consequently, for fixed $y>0$ 
$$(J^{(\v)}(t,(\cdot,y)):\ t>0)=(\sum_{j=1}^{\mathcal N(t)(x,y)}k_j:\ t>0)\overset{\text{\tiny dist.}}{=}(S^{(t)}:\ t>0).\ \ \CheckedBox$$
\

 \subsection*{Radon measures}
\

 \
 
 Let $X$ be a locally compact, Polish space.  We'll call  a measure $m:\B(X)\to [0,\infty]$  {\it Radon} 
  if $m(C)<\infty\ \forall\ C\subset X$ compact.
\

Let $\mu_n,\ m$ be Radon measures on  $X$. As in  \cite[Chapter 10]{Breiman}, we say that $\mu_n$ {\it converges weakly} to $m$
and write $\mu_n\xrightarrow[n\to\infty]{w}\ m$ if
$$\mu_n(U)\xrightarrow[n\to\infty]{}\ m(U)\ \forall\ U\in C_c(X)$$
where $C_c(X)$ denotes  the collection of continuous, $\Bbb R$-valued, 
compactly supported functions on $X$,  and where 
$\mu(U):=\int_XUd\mu$ for $U:X\to\Bbb R$ $\mu$-integrable.
This  abuse of notation will be useful when, in a dearth of clopen sets, we need to replace $1_A$  by $U\in C_c(X)$ with 
$\mu(A)=\mu(1_A)\approx\mu(U)$.
\

We'll call set $A\in\B(X)$ {\it almost open} if $m(\bdy A)=0$ and a function $f:X\to M$ (a metric space)
 {\it Riemann integrable} ({\tt RI}) if
  $$m(\{y\in X:\ f\ \text{not continuous at}\ y\})=0.$$

We'll denote the collections of $\text{almost open sets}$ and precompact, $\text{almost open sets}$ by 
 $\mathcal J(m)$ and $\mathcal J_c(m)$ respectively.
 \
 
As is well known $\mu_n\xrightarrow[n\to\infty]{w}\ m$ iff
$$\mu_n(A)\xrightarrow[n\to\infty]{}\ m(A)\ \forall\ A\in \mathcal J_c(m).$$ 

 \

Let $X$ be a  locally compact, Polish space. 
\

A collection $\mathcal H\subset C_c(X)$ is called
{\it measure determining} if for $\a,\ \b:\B(X)\to [0,\infty]$ Radon measures on $X$:
\par $\a(C)=\b(C)\ \forall\ C\in\mathcal H$ implies $\a\equiv\b$. 
\

For example, if $\mathcal H\subset C_c(X)$ is closed under multiplication, separates points of $X$ and
\f $\forall\ U\subset X$ open, precompact, $\exists\ f\in\mathcal H,\ f\ge 1_U$, then by the Stone-Weierstrass theorem,
$\mathcal H$ is  measure determining.

 \proclaim{Weak Convergence Lemma}\ \ Let $(Z,\mu)$  be a  locally compact measure space. 
\

Suppose that $\mathcal H\subset C_c(X)$ is a  measure determining collection.
\

 If $\mu_n\ (n\ge 1)\ \&\ \mu$ are Radon measures on $X$ and
 $$\mu_n(A)\xrightarrow[n\to\infty]{}\mu(A)\ \forall\ A\in\mathcal H,$$
 then
 $\mu_n\xrightarrow[n\to\infty]{w}\mu.$\endproclaim\demo{Proof}
Fix a sequence  $(H_k:\ k\ge 1)$ of compact subsets of $X$ so that $H_k\subset H_{k+1}^o\uparrow Z$.
By compactness, $\forall\ U\in C_c(X),\ \exists\ k\ge 1$ so that $[U\ne 0]\subset H_k$.
\

\ \  Using Helly's theorem and diagonalization, we obtain $t_j\to\infty$ and a Radon measure $\frak m$ on $X$  so that
$\forall\ k\ge 1$,
$$\mu_{t_j}(U)\xrightarrow[j\to\infty]{}\frak m(U)\ \forall\ U\in C_c(X),\ [U\ne 0]\subset\ H_k,$$
and in particular $\forall U\in\mathcal H$.
\

By assumption,
$$\mu_{t_j}(A)\xrightarrow[j\to\infty]{}\mu(A)\ \forall\ A\in\mathcal H.$$
Since $\mathcal H$ is  measure determining, $\frak m\equiv \mu$. 
\

Thus,  the only accumulation point of $\{\mu_{t}:\ t\ge 1\}$ is $\mu$
and
$\mu_{t}\xrightarrow[t\to\infty]{w}\mu.$.\ \ \CheckedBox
 \subsection*{Riemann integrable transformations of locally compact measure spaces}
\

 A $\s$-finite, Polish measure space $(X,m)$ is {\it locally compact} if $X$ is locally compact and   $m:\B(X)\to [0,\infty]$ is so that  $m(C)<\infty\ \forall\ C\subset X$ compact.
 \
 
 \
 
 Let $(Y,\nu)\ \&\ (Y',\nu')$ be locally compact, separable measure spaces.
 \
 
It is standard that if $\pi:Y\to Y'$ is nonsingular in the sense that for $A\in\B(Y'),\ \nu(\pi^{-1}A)=0\ \iff\ \nu'(A)=0$, then 
\bul $\pi$ is {\tt RI} if and only if  $\pi^{-1}\mathcal J(\nu)\subset\mathcal J(\nu)$.

Consequently, the composition of nonsingular, {\tt RI} maps is {\tt RI}.
 \

 \subsection*{{\tt RI} semiflows}
 Let $(Y,\nu)$ be a  locally compact, measure space. A semiflow $\Phi:\Bbb R_+\to\text{\tt MPT}\,(Y,\nu)$
 is a {\it {\tt RI} semiflow} if $\forall\ R>0$,
 $(t,x)\mapsto \Phi_t(x)$ is a {\tt RI} map $(0,R)\x Y\to Y$ with respect to the measure $\text{\tt Leb}\x\nu$.
 \
 
 It is not hard to show that if $(Y,\nu,\Phi)$ is a {\tt RI} semiflow and
 $G:\Bbb R_+\x Y\to\Bbb G$  is a {\tt RI}, $\Phi$-cocycle where $\Bbb G\le\Bbb R^\kappa$ is a 
 closed subgroup of full dimension, then the skew product semiflow $(Z,\mu,\Phi^{(G)})$ is also {\tt RI}.
 \
 
Let  $(Z,\mu,\Psi)\ \&\ (Z',\mu',\Psi')$  be  {\tt RI} semiflows. 
We'll say that $(Z,\mu,\Psi)$ is a {\it {\tt RI} extension} of $(Z',\mu',\Psi')$ if there is a measure preserving, {\tt RI} {\it factor map}
 $\pi:Z\to Z'$ with $\pi\circ\Psi_t=\Psi_t'\circ\pi$ and that $(Z,\mu,\Psi)\ \&\ (Z',\mu',\Psi')$ are
 {\it {\tt RI}-conjugate} if there is an invertible  {\it {\tt RI}-factor map} $\pi:Z\to Z'$ with $\pi^{-1}$ {\tt RI}.
 \
 
 \subsection*{$(\mathcal S,b)$-mixing on a locally compact measure space}
 \

Local limit mixing  properties are defined for  {\tt RI} semiflows on locally compact measure  spaces. 
\

Let $0<p\le 2$, let $\mathcal S$ be  a symmetric, $p$-stable, globally supported, random variable  on 
 $\Bbb R^{\kappa}$ and let $b(t)\to\infty$ be $\frac1p$-regularly varying.
 \
 
 We call the {\tt RI}  skew product semiflow $(Z,\mu,\Phi^{(G)})$  
 \sms -- (integrated) {\it  $(\mathcal S,b)$-mixing} if  
 $\forall\ A,\ B\in\mathcal J_c(\mu)$,
  \begin{align*}\tag{\Bat}\label{Bat}
 b(t)^\kappa\mu(B\cap\Phi_t^{(G)-1}Q_{k(t)}A)\xrightarrow[t\to\infty,\ k(t)\in\Bbb G,\ 
    \frac{k(t)}{b(t)}\to z]{}\ f_\mathcal S(z)\mu(A)\mu(B);
  \end{align*}

  \sms -- {\it conditionally  $(\mathcal S,b)$-mixing} if  for $ A\in\mathcal J_c(\mu)$,
 \begin{align*}\tag{\dsaeronautical}\label{dsaeronautical}
 b(t)^\kappa\widehat{\Phi_t^{(G)}}(1_{Q_{k(t)}A})(x)\xrightarrow[t\to\infty,\ k(t)\in\Bbb G,\ 
    \frac{k(t)}{b(t)}\to z]{}\ f_\mathcal S(z)\mu(A)\ \text{a.e.}.
  \end{align*} 
Here $Q_{k(t)}$ denotes a deck holonomy as defined on p. \pageref{deck}.
\

  In order to use weak convergence of Radon measures, we'll need following  standard (equivalent) reformulations of (\Bat) and 
  (\dsaeronautical)  in terms of $C_c(Z)$.
  For $U,\ V\in C_c(Z)$,
 {\small\begin{align*}&\tag{\Bat$_{C}$}
 b(t)^\kappa\mu(U\cdot V\circ Q_{k(t)}\circ\Phi_t^{(G)})\xrightarrow[t\to\infty,\ k(t)\in\Bbb G,\ 
    \frac{k(t)}{b(t)}\to z]{}\ f_\mathcal S(z)\mu(U)\mu(V);\\ &\tag{\dsaeronautical$_{C}$}
 b(t)^\kappa\widehat{\Phi_t^{(G)}}(U\circ Q_{k(t)})(x)\xrightarrow[t\to\infty,\ k(t)\in\Bbb G,\ 
    \frac{k(t)}{b(t)}\to z]{}\ f_\mathcal S(z)\mu(U)\ \ \text{a.e.}.
  \end{align*}}

Let $\mathcal H\subset C_c(Z)$.
We'll say that  the conditional   
$(\mathcal S,b)$-mixing  of the {\tt RI}  skew product semiflow $(Z,\mu,\Phi^{(G)})$ is {\it $\mathcal H$-uniform}   
if   for $U\in\mathcal H$,
$M>0\ \&\ E\subset\ Z$ compact, we have
$$\sup_{x\in E}|b(t)^\kappa\widehat{\Phi_t^{(G)}}(U\circ Q_{k(t)})(x)- f_\mathcal S(z)\mu(U)|\xrightarrow[t\to\infty,\ k(t)\in\Bbb G,\ \|k(t)\|\le Mb(t),\ 
    \frac{k(t)}{b(t)}\to z]{}\ 0.$$

 Theorem 1 (below) gives sufficient conditions for $(\mathcal S,b)$-mixing (both integrated and conditional).
 
  \subsection*{Remarks}
  \
  
  (o) Both conditional and integrated local limit mixing are invariant under {\tt RI} conjugacy.
  \
  
(i) By taking $k(t):=0$ in (\Bat) we see that  $(\mathcal S,b)$-mixing implies {\it Krickeberg mixing}
as in \cite{Krickeberg-recent} ($\mathcal J_c(\mu)$-{\it mixing}  in \cite{Friedman}):
\begin{align*}
 b(t)^\kappa\mu(B\cap\Phi_t^{(G)-1}A)\xrightarrow[t\to\infty]{}\ \mu(A)\mu(B)\ \forall\ A,\ B\in\mathcal J_c(\mu).
  \end{align*}
\

(ii) The  local limit theorems ({\tt Int-LLT}) and  ({\tt Con-LLT}) 
for $G$ with respect to $\mathcal R=\mathcal J_c(\nu)$ (as on p. \pageref{Int-LLT}) 
are consequences of
(\Bat) and (\dsaeronautical) (respectively). 

 \
 
The mixing lemma  shows that the reverse implications also hold.

\

 \proclaim{Mixing Lemma}\ \ Let $(Z,\mu,\Phi^{(G)})$  be a {\tt RI}  skew product semiflow on a locally compact measure space. Let
 $0<p\le 2$, let $\mathcal S$ be  a symmetric, $p$-stable, globally supported, random variable  on 
 $\Bbb R^{\kappa}$ and let $b(t)\to\infty$ be $\frac1p$-regularly varying.
\

Suppose that $\mathcal H\subset C_c(X)$ is a  measure determining collection.

\sms {\rm (i)}\ If {\rm (\Bat)$_C$} holds $\forall\ A,\ B\in\mathcal H$, then $\Phi^{(G)}$ is $(\mathcal S,b)$-mixing.
\

\sms {\rm (ii)}\ If
{\rm (\dsaeronautical)$_C$} holds $\forall\ A\in\mathcal H$, then $\Phi^{(G)}$ is conditionally  $(\mathcal S,b)$-mixing.
\endproclaim
\demo{Proof}\ \ We'll only prove (ii), the proof of (i) being similar and simpler.
\

Fix $k:\Bbb R_+\to\Bbb G$ measurable, $\frac{k(t)}{b(t)}\to z$ as $t\to\infty$. We'll apply the weak convergence lemma to the 
{\it transition measures} $\mu^{(x)}_t$ defined (using standard disintegration theory) for a.e. $x\in X\ \&\ t>0$ by
$$\mu^{(x)}_t(U):=b(t)^\kappa\widehat{\Phi_t^{(G)}}(U\circ Q_{k(t)})(x) \ \ (U\in C_c(X)).$$
These are finite measures and the map $$(t,x)\in\Bbb R_+\x Z_0\mapsto\mu_t^{(x)}\in\{\text{Radon measures on $Z$}\},$$ is measurable.
\

By assumption, for a.e. $x\in Z$,
$$\mu^{(x)}_t(U)\xrightarrow[t\to\infty]{}f_\mathcal S(z)\mu(U)\ \forall\ U\in\mathcal H$$
whence by the weak convergence lemma\ \   $\mu^{(x)}_t\xrightarrow[t\to\infty]{w}f_\mathcal S(z)\mu$ \ \ \ 
\

and $(Z,\mu,\Phi^{(G)})$ is conditionally $(\mathcal S,b)$-mixing.\ \Checkedbox

   \subsection*{Smooth cocycles over suspended semiflows } 
\

Given $f:Y\to\Bbb R^\kappa$ measurable and {\it flow integrable} in the sense that $\int_0^{\frak r(x)}|f(x,t)|dt<\infty$ a.e., 
we can define (as in 
\cite{Waddington},\ \cite{Iwata} and \cite{DN}) the {\it smooth cocycle} 
$F^{(f)}:\Bbb R_+\x Y\to\Bbb G$  by
\begin{align*}
 F^{(f)}(t,(x,y)):=\int_0^tf\circ\Phi_s(x,y)ds.
\end{align*} This is  a $\Phi$-cocycle.
It is  cohomologous to the jump cocycle $J^{(\varphi)}$ where
$\varphi:X\to\Bbb R^\kappa$ is defined by
$\varphi(x):=\int_0^{\frak r(x)}f(x,y)dy$:
\begin{align*}F^{(f)}(t,(x,u))=J^{(\varphi)}(t,(x,u))-\mathcal E(x,u)+\mathcal E\circ\Phi_t(x,u) 
\end{align*}
where
$\mathcal E(x,u):=\int_0^uf(x,s)ds.$
\

In case the $\Phi\ \&\ f$ are {\tt RI}, the skew product semiflows $\Phi^{(F)}$ and $\Phi^{(J)}$
on $(Y\x\Bbb R^\kappa,\nu\x m_{\Bbb R^\kappa})$ are {\tt RI} conjugate and 
  share the same $(\mathcal S,b)$-mixing properties.

\subsection*{Conditional rational weak mixing}
  \
  
  We call the measure preserving transformation $(Z,\mu,R)$ {\it conditionally rationally weakly mixing}
  if there exist constants $u_n>0$ so that with $a(n):=\sum_{k=1}^nu_k$,
  \begin{align*}\tag{{\tt P}}
    \frac1{a(N)}\sum_{k=1}^N|&\widehat{R}^k1_C-u_k\mu(C)|\\ &\xrightarrow[N\to\infty]{} 0\ \text{a.e.}\ \forall\ C\in\B(Z),\ 
    0<\mu(C)<\infty.
  \end{align*}
 This property was considered in \cite{nice}.
\

\

Let $0<p\le 2$, let $\mathcal S$ be  a symmetric, $p$-stable, globally supported, random variable  on 
 $\Bbb R^{\kappa}$ and let $b(t)\to\infty$ be $\frac1p$-regularly varying so that 
 $$a(N):=\sum_{n= 1}^N\frac1{b(n)^\kappa}\xrightarrow[N\to\infty]{}\ \infty$$
 (i.e. $\kappa=1\ \&\ 1\le p\le 2$ or $\kappa=2\ \&\ p=2$).
 \
 
 We'll see that if $(Z,\mu,\Phi^{(G)})$  is an ergodic, $(\mathcal S,b)$-mixing, skew product semiflow then
 each $\Phi_t^{(G)}$ is conditionally rationally weakly mixing. Theorem 2 gives sufficient conditions for this which are more general
 than those of theorem 1.

\subsection*{Transfer operator of a suspended semiflow}
\

Let $(Y,\nu,\Phi)=(X,m,T)^{\frak r}$ be a suspended semiflow. For $t>0$, the transfer operator of $\Phi_t$ is the operator 
$\widehat{\Phi}_t:L^1(\nu)\to L^1(\nu)$ defined by
$$\int_A\widehat{\Phi}_tfd\nu=\int_{{\Phi}_t^{-1}A}fd\nu\ \ (f\in L^1(\nu),\ A\in\mathcal C).$$
For $f=1_{B\x I}$ with $B\x I\in\mathcal C$, it is  given by
\begin{align*}\widehat{\Phi}_t(1_{B\x I})(\om,y)&=\sum_{n=0}^\infty\T^n(1_{[\mathcal N(t)=n]\cap {B\x I}}(\cdot,y))(\om)\\ &
=\sum_{n=0}^\infty\T^n(1_{B\cap [\frak r_n(\cdot)\in I+t-y]})(\om). 
\end{align*}

\

\

\subsection*{Fibered systems}
A {\it fibered system} is a quadruple $(X,m,T,\alpha)$ where $(X,m,T)$ is a non-singular transformation 
 and $\alpha$ is a countable, measurable partition which generates $\B$ under $T$ 
  ($\s(\{T^{-n}\a:\ n\ge 0\})=\B$)  such 
that  $T|_a$ invertible and nonsingular for $a\in\alpha$.
\

If $(X,m,T,\alpha)$ is a fibered system, then for $n\ge 1$, so are \f $(X,m,T,\alpha_n)\ \&\ (X,m,T^n,\alpha_n)$ where $\alpha_n=\a_0^{n-1}:=\bigvee_{k=0}^{n-1}T^{-k}\alpha$.

\

The collection of {\it $\a$-cylinders} of $(X,m,T,\alpha)$ is
$$\label{cylinders}\mathcal C_\a:=\bigcup_{n\ge 1}\a_n.$$
\

For $n\ge 1$, there are $m$-nonsingular inverse branches of $T^n$
denoted $v_a:T^na\to a\ (a\in\a_n)$ with Radon Nikodym derivatives
$$v_a':={dm\circ v_a\over dm}:T^na\to\Bbb R_+.$$
\

Fibered systems $(X,m,T,\alpha)$ and $(X',m',T',\alpha')$ are {\it isomorphic} if there is a measure theoretic isomorphism
$\pi:(X,m,T)\to (X',m',T')$ so that $\pi(\a)=\a'$. 
\

Any fibered system $(X,m,T,\alpha)$ is isomorphic to a {\it subshift} fibered system, that is one where 
$X\subset S^\Bbb N$  is $T$-invariant and closed  with respect to the  Polish topology on $S^\Bbb N$ (where  $S\subset\Bbb N$);  $T:S^\Bbb N\to S^\Bbb N$  is
the shift and $\a=\{[s]=\{(x_1,x_2,\dots)\in X:\ x_1=s\}:\ s\in S\}$.%
\

Thus, in a subshift fibered system $(X,m,T,\alpha)$, $X$ is a Polish space and $T:X\to X$ is continuous with $\a$  a collection of clopen sets in $X$. 
\

If $S$ is finite, then $X$ is compact, but if not, $X$ may not even be locally compact. 
\subsection*{Locally compact  fibered systems}
\

We'll need these in order to  consider  $(\mathcal S,b)$-mixing properties
of their skew product extensions and suspensions.
\

The fibered  system $(X,m,T,\alpha)$ is {\it locally compact} if $(X,m,T)$ is locally compact, $T$ is {\tt RI} and
$\C_\a\subset\mathcal J_c(m)$  where $\C_\a$ denotes the collection of $\a$-cylinders as on p. \pageref{cylinders}.
\

The fibered system $(X,m,T,\alpha)$ is {\it Markov} if $Ta\in\s(\a)\ \forall\ a\in\a$. In this case, so are
 $(X,m,T,\alpha_n)\ \&\ (X,m,T^n,\alpha_n)$.
\

An {\it interval map} is a fibered system $(X,m,T,\alpha)$ with  $X\subset\Bbb R$ a bounded interval, $m$  Lebesgue measure and $\alpha$  a 
partition $\mod m$ of $X$ 
into open intervals so that for each $a\in\alpha$, $T|_a$ is (the restriction of) a bi-absolutely continuous homeomorphism.
It is locally compact.
\

Next,  Markov fibered systems have locally compact representations:
\proclaim{Interval map  lemma}\ \ Let $(X,m,T,\alpha)$ be a Markov fibered system, then $(X,m,T,\alpha_2)$
is isomorphic to a Markov interval map.\endproclaim 
\demo{Proof}\ \ WLOG 
$(X,m,T,\alpha)$ is   a  subshift fibered system as above with
$X\subset S^\Bbb N$   $T$-invariant and closed and so that $s,t\in S\subset\Bbb N,\ s<t\ \implies\ m([s])\ge m([t])$.

\

Let $Y:=[0,1]\ \&\ \mu:=\text{\tt Leb}$.
\

Construct a partition  $\b=\{B_s:\ s\in S\}$  of $Y$ into open intervals
so that

\sms (i) $t,\ u\in S,\ t<u\ \implies\ B_t<B_u$ (i.e. $x<y\ \forall\ x\in B_t,\ y \in B_u$)
\sms (ii) $\mu(B_s)=m([s])\ \forall\ s\in S$;

Continue, constructing  partitions $\b_n\ (n\ge 1)$ of Y into open intervals of form
$$\b_n:=\{B_{(k_1,k_2,\dots,k_n)}:\ \u k=(k_1,k_2,\dots,k_n)\in S^n\}$$ satisfying
\sms (a) $\u k\prec\u\ell$ (lexicographically) $\Rightarrow\ B_{\u k}<B_{\u\ell}$;
\sms (b) $\mu(B_{\u k})=m([\u k])$.
\sms (c)  $\forall\ n\ge 1,\ \u k=(k_1,k_2,\dots,k_n)\in S^n$,
$$B_{\u k}=\bigcupdot_{i\ge 1}B_{\u k,i}.$$
Define a map $\pi:X\to [0,1]$ by
$$\bigcap_{n\ge 1}\overline{B_{(x_1,x_2,\dots,x_n)}}=\{\pi(x)\}.$$
The map is continuous and $m\circ\pi^{-1}=\mu$. 
It is injective on $\pi^{-1}(Y_0)$ with
$$Y_0:=Y\setminus\bigcup_{n\ge 1,\ b\in\b_n}\bdy b$$
-- a dense $G_\d$ set of full measure and, indeed, $\pi:\pi^{-1}(Y_0)\to Y_0$ is an homeomorphism.
\

Thus  $\tau:Y_0\to Y_0$ defined by
$\tau(y):=\pi(T(\pi^{-1}(y))$ is continuous and $\pi:(X,m,T,\alpha_2)\to (Y,\mu,\tau,\b_2)$ is an isomorphism of fibered systems
\

It remains to show that $(Y,\mu,\tau,\b_2)$ is an interval map. 
\

To this end, 
\Par\ {\it For $k,\ell\in S$ so that $m([k,\ell])=\mu(B_{k,\ell})>0$,
$\tau:B_{k,\ell}\cap Y_0\to B_\ell\cap Y_0$ extends to an increasing, continuous map $\tau_{k,\ell}:B_{k,\ell}\to B_\ell$.}
\pf \ For such $k,\ell\in S$, by the Markov property, $T:[k,\ell]\to [\ell]$ is surjective, whence also
  $\tau:B_{k,\ell}\cap Y_0\to B_\ell\cap Y_0$.
  \
  
  First,  $\tau:B_{k,\ell}\cap Y_0\to B_\ell\cap Y_0$ is strictly increasing because:
\par if $x,y\in B_{k,\ell}\cap Y_0,\ x<y$, then $\exists\ n\ge 1$ and
  $\u x,\ \u y\in S^n$ with $x_j=y_j\ \forall\ 1\le j<n\ \&\ x_n<y_n$ so that
  \par $x\in B_{k,\ell,\u x}\ \&\ y\in B_{k,\ell,\u y},$ whence $\tau(x)\in B_{\ell,\u x}\ \&\ \tau(y)\in B_{\ell,\u y}$
 and $\tau(x)<\tau(y)$.
  \
  
  Next, for $z\in B_{k,\ell}$, set
  $$\tau_{k,\ell}(z):=\lim_{x\to z-,\ x\in B_{k,\ell}\cap Y_0}\tau(x).$$
  
We claim that $\tau_{k,\ell}:B_{k,\ell}\to B_\ell$ is continuous and increasing.

  \f{\tt Proof of increasing}:
  \begin{align*}
x,y\in B_{k,\ell},\ x<y\ &\Rightarrow\ \exists\ u,v\in (x,y)\cap Y_0,\ u<v\\ & \implies\ \tau_{k,\ell}(x)\le \tau_{k,\ell}(u)<\tau_{k,\ell}(v)\le \tau_{k,\ell}(y).
  \end{align*}
  \
  
   \f{\tt Proof of continuity}:
\ Being monotone, either $\tau_{k,\ell}:B_{k,\ell}\to B_\ell$  is continuous or there is a nonempty, open interval 
$J\subset B_\ell$ with   $\tau_{k,\ell}(B_{k,\ell})\subset B_\ell\setminus J$  contradicting 
  $\tau_{k,\ell}(B_{k,\ell})\supset B_\ell\cap Y_0$.\ \ \Checkedbox\ \P
  
  \
  
  It is now standard that $(Y,\mu,\tau,\b_2)$ is a Markov interval map.\ \ \Checkedbox

\subsection*{Quasicompact action}\label{quasicompact}
\

Let $L\subset L^\infty$ be a Banach space  which is 
{\it admissible} in the sense  that $\|\cdot\|_L\ge\|\cdot\|_\infty$. We say that {\it $\T$ acts  quasicompactly on $L$} if $\T L\subset L$ and $\exists\ M\ge 1,\ \rho\in (0,1)$ so that
$$\|\T^n f-\int_Xfdm\|_L\le M\rho^n\|f\|_L\ \ \forall\ n\ge 1,\ f\in L.$$
This definition of {\tt quasicompact action} differs from the one in \cite{HH}, but coincides with it when the fibered system is weakly mixing.
\

The fibered systems $(X,m,T,\a)$ have transfer operators which act quasicompactly on an associated admissible Banach space denoted $L_T$.
\subsection*{Example 1:\ {\tt AFU} maps}
\

An {\it {\tt AFU} map} is an  interval map   $(X,m,T,\alpha)$ where  for each $a\in\alpha$, $T|_a$ is (the restriction of) a  $C^2$ diffeomorphism 
$T:\overline{a}\to \overline{Ta}$
 satisfying in addition:
                
\begin{align*}&\tag{A}\sup_X\tfrac{|T''|}{(T')^2}<\infty,\\ &
 \tag{F}T\a:=\{TA:A\in\alpha\}\ \text{\rm is finite},\\ &\tag{U}\inf_{x\in a\in\alpha}|T'(x)|>1.
\end{align*}
An {\tt AFU} map with finite $\alpha$ is known as a {\it Lasota-Yorke} map after \cite{LY}. 
\

It is known that

\bul \cite{Rychlik},\ \cite{Roland}:  an {\tt AFU} map has a Lebesgue-equivalent, invariant probability whose density has bounded variation
and which is uniformly bounded below;
\bul \cite{Rychlik}: for $T$ a mixing, {\tt AFU} map,  $\T$ acts  quasicompactly on  $L_T=\text{\tt BV}$, the space of functions on the interval $X$ with bounded variation; 
this quasicomptacness persists when $m$ is replaced by the absolutely continuous, $T$-invariant probability;
\bul  \cite{Rychlik}: a  topologically mixing {\tt AFU} map is exact and (\cite{coeffts})
the generated stochastic process $(\alpha\circ T^n:\ n\ge 1)$ is  {\it reverse $\phi$-mixing} in the sense that
\begin{align*}
\label{fatso}\phi_-(n):=
\sup_{k\ge 1}\{m(A|B)-m(A)|:\ A\in\a_0^{k-1},\ B\in\a_{k+n}^\infty\}\xrightarrow[n\to\infty]{}\ 0
\end{align*}
 where  $\a_J^{K}:=\s(\{\a\circ T^j:\ J\le j<K+1\})$ for $0\le J<K\le\infty$. 
 \
 
 Indeed, here, $\phi_-(n)\xrightarrow[n\to\infty]{}\ 0$ 
 exponentially.

\subsection*{Example 2:\ Gibbs-Markov maps}
\

A {\it Gibbs-Markov {\rm ({\tt G-M})} map} $(X,m,T,\alpha)$ is a fibered system 
 which is {\it Markov} in the sense that 
 \begin{align*}
   \tag{a}Ta\in\s(\alpha)\ \mod m\ \ \forall\ a\in\a
 \end{align*}

and which satisfies 
 \begin{align*}
\tag{b}\inf_{a\in\alpha}m(Ta)>0 \end{align*}
 and, for some $\th\in (0,1)$
 \begin{align*}\tag{$\frak g_\th$}\sup_{n\ge 1,\ a\in\alpha_n}D_{T^na,\th}(\log v_a')<\infty
\end{align*}
where, for $A\subseteq\ X,\ H:A\to\Bbb R$
 \begin{align*}
 D_{A,\th}(H):=\sup_{x,y\in A}\frac{|H(x)-H(y)|}{\th^{t(x,y)}}<\infty
 \end{align*}

with $t(x,y)=\min\{n\ge 1:\ \a_n(x)\ne \a_n(y)\}\le\infty$.

A {\tt G-M} map on a finite state space is a subshift of finite type ({\tt SFT}).
\

It is known (see \cite{AD}) that
\bul a {\tt G-M} map has a $m$-equivalent, invariant probability whose density $h$ is $\th$-{\it H\"older continuous} on $X$
in the sense that $D_{X,\th}(h)<\infty$ where, for $A\subseteq\ X$
 \begin{align*}
 D_{A,\th}(h):=\sup_{x,y\in A}\frac{|h(x)-h(y)|}{\th^{t(x,y)}}<\infty
 \end{align*}

where $t(x,y)=\min\{n\ge 1:\ x_n\ne y_n\}\le\infty$;
\bul  for $T$ a mixing, {\tt G-M} map,  
$\T$ acts quasicompactly on $L_T=H_\th$, the space of 
$\th$-H\"older continuous functions $X\to\Bbb R$ with
$\th$ as in ($\frak g_\th$) equipped with the norm 
$$\|f\|_{H_\th}:=\|f\|_\infty+D_{X,\th}(f).$$
\bul if the invariant density is bounded below, the quasicompactness persists when $m$ is replaced by the absolutely continuous, $T$-invariant probability;
\bul  a  topologically mixing {\tt G-M} map is exact and the stochastic process $(\alpha\circ T^n:\ n\ge 1)$ is exponentially 
{\tt  continued fraction mixing},
\bul a Markov {\tt AFU} map is also a {\tt G-M} map.

\subsection*{Aperiodicity and non-arithmeticity}
\

Let $\Bbb H$ be a locally compact, Abelian, Polish group.
\

We call the measurable $G:X\ \to\Bbb H$ :
\bul {\it non-arithmetic} if
$\nexists$ solution to
\begin{align*}&G=k+g-g\circ T,\ \ g:X\ \to\Bbb H\ \ \text{measurable},\\ & k:X\ \to \Bbb K\ \ 
\text{measurable, where $\Bbb K$ is a proper subgroup of $\Bbb H$}. 
\end{align*}
and
\bul {\it aperiodic} if
$\nexists$ solution to
\begin{align*}&G=k+g-g\circ T,\ \ g:X\ \to\Bbb H\ \ \text{measurable},\\ & k:X\ \to \Bbb K\ \ 
\text{measurable, where $\Bbb K$ is a proper coset of $\Bbb H$}. 
\end{align*}
\

\subsection*{Geodesic flows on cyclic covers }\label{geodesic}
\

We denote by $(U(M),\L,g)$, the {\tt geodesic flow} on the {\tt unit tangent bundle} $U(M)$ of the  {\tt hyperbolic surface} $M$ equipped with
the $g$-invariant {\tt Liouville measure} $\L$  (for definitions see \cite{Hopf},\ \cite{Rees},\ \cite{nice}). 
\

The hyperbolic surface $V$ is a {\it cyclic cover} of the compact  hyperbolic surface $M$ if there is
a covering map $p:V\to M$ and a monomorphism $\g:\Bbb Z\to\text{\tt Isom}\,(V)$ (hyperbolic isometries of $V$), so that
for $y\in V,\ p^{-1}\{p(y)\}=\{\g(n)y:\ n\in\Bbb Z\}$. 
\

It was shown in \cite{nice} that such $(U(V),\L,g)$ is rationally weakly mixing, and rationally ergodic of order $2$ (as defined there). 
This is strengthened in
theorem 3 (below) where we establish $(\mathcal S,b)$-mixing (with\  $\mathcal S$ centered Gaussian and $b(t)\propto\sqrt n$) along with 
{ multiple Krickeberg mixing} and 
{ rational weak mixing  of order $2$}.
\section*{\S2 Results}

\subsection*{Set-up}
\

We'll prove local limit mixing theorems (Theorems 1 and 2 below) for $J^{(\varphi)}$ (defined above in (\Football) on p. \pageref{Football}) with 
respect to the semiflow $(Y,\nu,\Phi)=(X,T,\mu,\alpha)^\frak r$.
with 
$(X,T,\mu,\alpha)$ a topologically mixing, probability preserving fibered system,
$\frak r:X\to\Bbb R_+$ and $\varphi:X\ \to\Bbb G$ measurable satisfying certain conditions listed below.
\

\subsection*{General assumptions}
\

For all results, we assume that 
\bul $(X,m,T,\alpha)$ is probability preserving,  mixing,  locally compact and either a {\tt G-M} map or an {\tt AFU} map;
\bul $\frak r\in L_T$ (see p. \pageref{quasicompact})\ $\&\ \frak r>0$;
\bul\label{osc} $\mathcal D^{(T)}_\alpha(\varphi)<\infty$ where
\par\ \ \ \  for $(X,m,T,\alpha)$ a {\tt G-M} map, $$\mathcal D^{(T)}_\alpha(\varphi):=\sup_{a\in\alpha}\,D_{a,\th}(\varphi)<\infty$$
\par\ \ \ \  and for $(X,m,T,\alpha)$ an {\tt AFU} map, $$\mathcal D^{(T)}_\alpha(\varphi):=\sup_{a\in\alpha}\,\text{\tt Lip}_{a}(\varphi)<\infty.$$

Here, for $A\subset X$,
$$\text{\tt Lip}_A(\varphi):=\sup\,\{\frac{\|\varphi(x)-\varphi(y)\|}{|x-y|}:\ x,y\in A\}.$$
Note that for $A\subset X$ an interval,
$$\bigvee_A\v:=\sup_{n\ge 1,\ t_0<\dots< t_n,\ t_0,\dots, t_n\in A}\,\sum_{k=0}^{n-1}\|\v(t_{k+1})-\v(t_k)\|\le m(A)\text{\tt Lip}_{A}(\varphi).$$
Also, by (U), if $\vartheta\in (0,1)$ is sufficiently large, then 
$$D_{a,\vartheta}(\varphi)\le \text{\tt Lip}_{a}(\varphi)\ \forall\ a\in\alpha.$$
\

\subsection*{Weak Independence}
\

Our {\tt weak independence} assumption is that for some $\frak c>0$ we have
\begin{align*}\tag{{\scriptsize\faAmbulance}}\label{faAmbulance}
 m([\v_n\in U,\ \frak r_n\in J])&\le \frak cm([\v_n\in U])m([\frak r_n\in J])\\ & 
 \forall\ n\ge 1,\ U\subset\Bbb G,\ \ J\in\Bbb R_+\ \text{open subsets}.
\end{align*}

Under the general assumptions, 
a restricted version of ``weak 
independence'' holds whenever $\frac{\v_n}{b(n)}$ converges in distribution to a non-normal stable random variable.
See ({\scriptsize\faPuzzlePiece}) on p. \pageref{Puzzle}.
But we do not know if this suffices for ({\scriptsize\faAmbulance}).
\

On the other hand, this ``weak independence'' cannot hold for the whole {\it roof-} and {\it dispacement-}  processes ($(\frak r\circ T^n:\ n\ge 0)$ and 
$(\varphi\circ T^n:\ n\ge 0)$ respectively) without their independence. Indeed:
\Smi\ \  If,   for some $M>0$,
$$m(A\cap B)\le Mm(A)m(B)\ \forall\ A\in \s(\{\frak r\circ T^n:\ n\ge 0\}),\ 
B\in \s(\{\varphi\circ T^n:\ n\ge 0\}),$$ 
then 
$$m(A\cap B)=m(A)m(B)\ \forall\ A\in \s(\{\frak r\circ T^n:\ n\ge 0\}),\ 
B\in \s(\{\varphi\circ T^n:\ n\ge 0\}).$$

To see \smiley, write 
$\mathcal A_\frak r:=\s(\{\frak r\circ T^n:\ n\ge 0\})\ \&\ \mathcal A_\varphi:=\s(\{\varphi\circ T^n:\ n\ge 0\})$
define the probability $\mu\in\mathcal P(X,\mathcal A_\frak r\vee \mathcal A_\varphi)$ by
$\mu(A\cap B):=m(A)m(B).$
\

The mixing properties of 
$(X,m,T)$ ensure that 
$(X,\mathcal A_\frak r\vee \mathcal A_\varphi,\mu,T)$ is an {\tt EPPT} where $\mathcal A_\frak r\vee \mathcal A_\varphi:=\s(\mathcal A_\frak r\cup\mathcal A_\varphi)$).
\

So is $(X,\mathcal A_\frak r\vee \mathcal A_\varphi,m|_{\mathcal A_\frak r\vee \mathcal A_\varphi},T)$
(as a factor of $(X,m,T)$). 
\

Thus $\mu\ \&\ m|_{\mathcal A_\frak r\vee\mathcal A_\varphi}$ are either disjoint or identical; the former being eliminated by the assumption.\ \Checkedbox\ \smiley

\subsection*{Aperiodicity}
We'll mainly consider  the {\tt joint aperiodicity} assumption:
\sms ({\tt J-Ap})\label{J-Ap}  \ \ \ \ $(\frak r,\varphi):X\ \to\Bbb R\x\Bbb G$ is aperiodic

\

In the sequel, we'll prove the  semiflow {\tt LLT} under ({\tt J-Ap}).
\

\subsection*{Distributional assumptions}
\

For some $0<p\le 2$ and some $\frac1p$-regularly varying sequence $(b(n):\ n\ge 1)$,
\begin{align*} 
 \tag{\dsliterary}\frac{\varphi_n}{b(n)}\xrightarrow[n\to\infty]{\text{\tiny dist.}}\ \mathcal S
\end{align*}  where 
 $\mathcal S$ is  a globally supported, symmetric $p$-stable (S$p$S) $\rv$ on $\Bbb R^{\kappa}$.
 \
 
 We'll consider the following separate cases of (\dsliterary):
 \sms ({\tt CN}) the { classical normal case} where $E(\|\varphi\|^2)<\infty,\ b(n)=\sqrt n$ and $\mathcal S$ is normal,
\sms ({\tt NNS})  the { non-normal, stable case} where $p\ne 2$.
\

\f{\bf Remark}\ \ We consider it likely that our results in the {\tt NNS} case  extend to the case where the displacement limit random variable $\mathcal S$ is strictly stable.
 However, the methods in the present paper and their dependencies rely on the (stronger) symmetry assumption.

\subsection*{Statements}
\

\proclaim{Theorem 1: semiflow  local limit mixing}\ 
\

\

{\rm (i)}\   Under the general assumptions, {\tt (J-Ap)}, and {\tt (CN)} with $E(\|\v\|^4)<\infty$,
the skew product semiflow is conditionally
  $(\mathcal S,b)$-mixing and this is $\mathcal H$-uniform with 
  $$\mathcal H:=L_T\x C_c(\Bbb R_+)\x C_c(\Bbb G)\subset C_c(Z).$$
 \
 
{\rm (ii)}\   Under the general assumptions, {\tt (J-Ap)}, weak independence and {\tt (NNS)} or  {\tt (CN)} with $E(\|\v\|^4)=\infty$
the skew product semiflow is   $(\mathcal S,b)$-mixing.
\

\endproclaim
\subsection*{{\tt CTRW}  -- a  toy model for theorem 1 }
\ 

\proclaim{{\tt LLT} for {\tt CTRW}s}\ \ Let $(S^{(t)}:\ t>0)$ be a {\tt CTRW} on $\Bbb G\le\Bbb R^\kappa$, a subgroup of full dimension, with  jump  random variables $(X_k:\ k\ge 1)$ and intensity $\l>0$.

\

Suppose that $\frac{S_n}{b(n)}\xrightarrow[t\to\infty]{\text{\tiny dist.}}\ \mathcal S$ where $p\in (0,2]$, $\mathcal S$ is  
a globally supported, symmetric, $p$-stable $\rv$ on $\Bbb R^{\kappa}$ and
$b(t)$ is $\frac1p$-regularly varying as $t\to\infty$, 
then for $U\subset\Bbb G$ a 
compact neighborhood of $0$ with $m_\Bbb G(\bdy U)=0$
\begin{align*}
b(\lambda t)^{\kappa}P(S^{(t)}\in k(t)+U)\xrightarrow[t\to\infty,\ \frac{k(t)}{b(\l t)}\to z\in\Bbb R^{\kappa}]{}\ m_\Bbb G(U)f_\mathcal S(z) 
\end{align*}
uniformly in $z\in K$ a compact subset of $\Bbb R^{\kappa}$ where $f_\mathcal S$ is the probability density function of $\mathcal S$.\endproclaim
For simple continuous time random walk, this is theorem 2.5.6 of \cite{L-L}. The general case uses the same methods.
In view of the discussion on p. \pageref{CTRW-jump}, the stronger statement in  Theorem 1(ii) can be obtained.
  
\

Assuming  ``recurrence'' (but not ``weak independence''), we have 

\proclaim{Theorem 2}\ 
\

  Under the general assumptions,  {\tt (J-Ap)}, and {\tt (CN)}  or {\tt (NNS)}, with
$$a(n):=\sum_{k=1}^n\frac1{b(k)^{\kappa}}\xrightarrow[n\to\infty]{}\infty,$$
the skew product semiflow is conditionally, rationally weakly mixing. 
\endproclaim

\proclaim{Theorem 3}
\

 For $V$  a  cyclic cover of a compact, hyperbolic surface, and $(U(V),\L,g)$
its geodesic flow, then for $\mathcal S\in\text{\tt RV}(\Bbb R)$  centered Gaussian and $b(t)=c\sqrt t$ (some $c>0$), 
\sms {\rm (o)} $(U(V),\L,g)$ is {\rm  $(\mathcal S,b)$ mixing} {\tt i.e.} for $A,\ B\in\mathcal J_c(\L)$, $k:\Bbb R_+\to\Bbb Z$  so that 
$\tfrac{k(t)}{\sqrt t}\xrightarrow[t\to\infty]{}z$, 
\begin{align*}\tag{{\scriptsize\faBus}}
b(t)\L(A\cap g_{t}Q_{k(t)}B)\ \xrightarrow[t\to\infty]{}\ \L(A)\L(B)f_\mathcal S(z);\end{align*}
\sms {\rm (i)} \ \ $(U(V),\L,g)$ is {\rm multiply Krickeberg mixing} {\tt i.e.}

for $N\ge 1,\ A,\ A_1,\dots,A_N\in\mathcal J_c(\L)$,
\begin{align*}\tag{{\scriptsize\faShip}}
\L(A\cap\bigcap_{i=1}^{N}g_{t_i}A_i)\ \sim\ \L(A)\prod_{i=1}^{N}\tfrac{\L(A_i)}{b(\D t_i)}\end{align*}
as $\u t=\{0=t_0<t_1<\dots<t_N\},\ \D t_i:=t_i-t_{i-1}\to\infty$, ($1\le i\le N$).
\sms {\rm (ii)}\ \  $(U(V),\L,g)$ is
{\rm rationally weakly mixing of order $2$}  {\tt i.e.} 
\

for each $\tau>0$ fixed,                                                                               
\begin{align*}\tag{\Yinyang}\label{Yinyang}\frac1{\log N}\sum_{k=1}^N|&\L(A\cap g_{\tau k}B\cap g_{2\tau k}C)-\tfrac{\L(A)\L(B)\L(C)}{b(\tau k)^2}|\\ &
\xrightarrow[N\to\infty]{}\ 0
\ \ \forall\ A,B,C\subset\ U(M)\ \ \text{bounded, measurable}.
\end{align*}
\endproclaim
\subsection*{Methods}
\

The proofs of  theorems 1 $\&$ 2 use an operator {\tt LLT} ({\tt Op-LLT}) for 
$(\frak r,\varphi):X\ \to\Bbb R\x\Bbb G$ with respect to $T$. 
In the ({\tt CN}) case, this boils down to the classical multidimensional {\tt LLT} (as in \cite{G-H}) . 
See \cite{AD} and \cite{ADSZ} for the versions considered here.
\

The {\tt Op-LLT}  is proved here  in the ({\tt NNS}) case (below).
For more information about operator stability $\&$ limits, see \cite{Jurek-Mason} and for an {\tt Op-LLT} in the independent case, see \cite{Doney}.
\

The proof of theorem 1 also uses an extended {\tt LLT} estimate for the renewal cocycle based on proposition 4 in \cite{Pene09}.
\

Theorem 3 is an application of theorem 1(i).

\subsection*{Plan of the paper}\ \ 
 
 \
 
 The proofs of  theorems 1 and 2 both use a splitting of the LHS of ({\tt Con-LLT}) (see p. \pageref{Con-LLT}) according to the value of the renewal cocycle.
 See  ({\scriptsize\faChainBroken}) on p. \pageref{faChainBroken}. 
 \

 The term I in ({\scriptsize\faChainBroken}) where the renewal cocycle takes moderate values is estimated using the the {\tt Op-LLT} in \S3.
 The remaining term  II is also treated in \S3  using a deviation lemma which is proved in the appendix (\S6) by means of an extended {\tt LLT} estimate also established
 there. 
 \
 
In \S4, we  prove the {\tt Op-LLT}, 
given an infinite divisibility lemma (p. \pageref{inf-div}), which in turn is established in \S5. 
\

\section*{\S3  Proofs}
\

\subsection*{Notations}\ \ In the sequel, we'll use the following notations:
\bul For sequences $a,\ b:\Bbb N\to\Bbb R$,\ $a(n)\approx b(n)$ means $a(n)-b(n)\xrightarrow[n\to\infty]{}\ 0$.
\bul For positive sequences $a,\ b:\Bbb N\to\Bbb R$,
\sms $a(n)\ll b(n)$ means that $\exists\ M>0$ so that $a(n)\le b(n)\ \forall\ n\ge 1$ large;
\sms $a(n)\gg b(n)$ means that  $b(n)\ll a(n)$; 
\sms $a(n)\asymp b(n)$ means that  $a(n)\ll b(n)$ and  $a(n)\gg b(n)$;
\sms $a(n)\sim b(n)$ means that $\tfrac{a(n)}{b(n)} \xrightarrow[n\to\infty]{}\ 1$;
\sms $a(n)\propto b(n)$ means that $\exists\ \lim_{n\to\infty}\tfrac{a(n)}{b(n)}\in\Bbb R_+$.
\bul For a sequence of $\Bbb R^d$-valued, measurable functions $(X_n:\ n\ge 0)$ defined on the probability space $(\Om,\mathcal F,P)$:
\sms $X_n\xrightarrow[n\to\infty]{\text{\tiny dist.}}\ Y$ (where $Y$ is a random variable on $\Bbb R^d$) means
$$\int_\Om g(X_n)dP \xrightarrow[n\to\infty]{}\ \Bbb E(g(Y))\ \forall\ g:\Bbb R^d\to\Bbb R,\ \text{bounded, continuous}.$$

random variables
\subsection*{Proof of theorem 1}
We consider the skew product semiflow
$(Z,\mu,\Phi^{(J)})$ over the suspended semiflow $(Y,\nu,\Phi)\cong\ (X,m,T,\a)^\frak r$ with 
$J=J^{(\v)},\ \v:X\to\Bbb G\ \&\ \frak r:X\to\Bbb R_+$. We'll use (\faLinux) on p. \pageref{faLinux}:
$$(Z,\mu,\Phi^{(J)})\cong\ (X\x\Bbb G,m\x m_\Bbb G,T_\v)^{\overline{\frak r}}$$ with $\overline{\frak r}(x,z):=\frak r(x)$.
\

It suffices to establish convergence for  the measure determining collection
 $$\mathcal H:=\C_\a\x\mathcal J_c(m_\Bbb R)\x\mathcal J_c(m_\Bbb G)\subset\mathcal J_c(\mu)$$
 where $\C_\a$ denotes the collection of $\a$-cylinders as on p. \pageref{cylinders}.

Fix $x:\Bbb R_+\to\Bbb G$ so that $\frac{x(t)}{b(\l t)}\xrightarrow[t\to\infty]{}\ z\in\Bbb R^{\kappa}$ where
$$\l:=\frac1{E(\frak r)}\in\Bbb R_+$$ and let 
$A\x I\x U\in\mathcal H\label{flowdisp}$.
\

We have
\begin{align*}\widehat{\Phi^{(J)}}_t(&1_{A\x I\x ( x(t)+U)})(x,z,y) =
\widehat{\Phi}_t(1_{(A\x I)\cap [-J^{(\varphi)}(t,\cdot)\in  x(t)-z+U]})(\om,y)\\ &=
\widehat{\Phi}_t(1_{(A\x I)\cap [J^{(\varphi)}(t,\cdot)\in  z-x(t)+V]})(\om,y)
\end{align*}
where $V:=-U$.
\

For each $M>0$, we have (as in \cite{nice}) the following splitting:

\begin{align*}\tag{{\scriptsize\faChainBroken}}\label{faChainBroken}\widehat{\Phi}_t(&1_{(A\x I)\cap [J^{(\varphi)}(t,\cdot)\in  z-x(t)+V]})(\om,y)\\ &
=\widehat{\Phi}_t(1_{(A\x I)\cap [J^{(\varphi)}(t,\cdot)\in  z-x(t)+V]}1_{[|\mathcal N(t)-\l t|\le M\sqrt t]})(\om,y)+\\ &
\ \ \ \ \ \ \ \ \ +\widehat{\Phi}_t(1_{(A\x I)\cap [J^{(\varphi)}(t,\cdot)\in  z-x(t)+V]}1_{[|\mathcal N(t)-\l t|>M\sqrt t]})(\om,y)\\ &
=:I(M,t,\om)+II(M,t,\om).
 \end{align*}

 The rest of the proof consists of establishing the following claims (I) and (II):
 \begin{align*}\tag{I}\label{I}\lim_{M\to\infty}\lim_{t\to\infty}b(\l t)^{\kappa} I(M,t,\om)=
 f_{\mathcal S}(z)\mu(A)\text{\tt Leb}\,(I)m_\Bbb G(U).
\end{align*}
uniformly on compact subsets of $Z$; which holds in all cases and where in  case (i) ({\tt CN}), we write $b(t):=\sqrt{t}$;
\sms \ \ \ and
\begin{align*}\tag{II}\label{II}&\varlimsup_{t\to\infty}b(t)^{\kappa}II(M,t,\cdot)\xrightarrow[M\to\infty]{}0\ \ \text{in case (i)};\\ &
\varlimsup_{t\to\infty}b(t)^{\kappa}E(II(M,t,\cdot))\xrightarrow[M\to\infty]{}0\ \ \text{in the cases (ii)}
\end{align*}
where in cases (o) $\&$ (i), the convergence is uniform on compact subsets of $Z$.
\

 Using  $(Y,\nu,\Phi)\cong\ (X,m,T,\a)^\frak r$ we have for $K\subset\Bbb N$,
 \begin{align*}\widehat{\Phi}_t(&1_{(A\x I)\cap [J^{(\varphi)}(t,\cdot)\in  z-x(t)+V]}1_{[\mathcal N(t)\in K]})(\om,y)\\ &
  =\sum_{n\in K}\widehat{\Phi}_t(1_{(A\x I)\cap [J^{(\varphi)}(t,\cdot)\in  z-x(t)+V]}1_{[\mathcal N(t)=n]})(\om,y)\\ &=
  \sum_{n\in K}\T^n(1_{A\cap [\varphi_n\in  z-x(t)+V,\ \frak r_n\in I+t-y]})(\om)
 \end{align*}
 whence
 \begin{align*}&I(M,t,\om)=\sum_{n\in\Bbb N,\ n= \l t\pm M\sqrt t}
 \T^n(1_{A\cap [\varphi_n\in  z-x(t)+V,\ \frak r_n\in I+t-y]})(\om)\\ &
 II(M,t,\om)=\sum_{n\in\Bbb N,\ |n-\l t|>M\sqrt t}
 \T^n(1_{A\cap [\varphi_n\in  z-x(t)+V,\ \frak r_n\in I+t-y]})(\om).
 \end{align*}

\

To prove (I), the independence assumptions are not needed. We shall need:
\bul in the ({\tt CN}) case, the local limit theorem for aperiodic $(\frak r,\varphi)$ as in \cite{G-H},

and
\bul in the ({\tt NNS}) case, an operator local limit theorem  ({\tt Op-LLT}) for  aperiodic  $(\frak r,\varphi)$

which we state now and prove in the sequel:

\

\subsection*{The operator local limit theorem statement}

\

\proclaim{Operator local limit theorem}
\

\ Suppose {\tt (J-Ap)}, and {(\tt NNS)} with $0<p<2$, 
 then there is a $\frac1p$-regularly varying sequence $(b(n):\ n\ge 1)$, a random variable $Z=(\mathcal S,\mathcal N)$
on $\Bbb R^{\kappa}\x\Bbb R$ with  a positive, continuous, probability density function, 
where $\mathcal S$ is  a globally supported, symmetric $p$-stable  $\rv$ on $\Bbb R^{\kappa}$, $\mathcal N$ is 
centered, Gaussian on $\Bbb R$, and independent of $\mathcal S$ so that for $A\in\mathcal C_\a$  and $I\subset\Bbb R$ an interval so 
that $A\x\ I\subset Y,\ \&\ U\subset\Bbb G$ a 
compact neighborhood of $0$ with $m_\Bbb G(\bdy U)=0$, 
\begin{align*}\tag{{\tt Op-LLT}}\label{Op-LLT}
\sqrt n b(n)^{\kappa}\T^n(&1_{A\cap [\varphi_n\in z_n+U,\ \frak r_n\in I+n{E(\frak r)}+r_n]})\\ & \underset{n\to\infty}{\text{\Large $\approx$}}\ 
f_Z(\zeta_n,\rho_n)m(A)|I|m_\Bbb G(U).
\end{align*} 
where $(z_n,r_n)\in\Bbb G\x\Bbb R$ with $z_n=b(n)\cdot\zeta_n,\ \ r_n=\sqrt n\cdot\rho_n$ with $(\zeta_n,\rho_n)=O(1)$.
\

For fixed $(z_n,r_n)\in\Bbb G\x\Bbb R$  as above, the convergence is uniform on $X$.\endproclaim
In the sequel, we'll refer to ({\tt Op-LLT}) as {\it the {\tt Op-LLT} corresponding to the distributional convergence}
 $$(\tfrac{\frak r_n-n{E(\frak r)}}{\sqrt n},\tfrac{\varphi_n}{b(n)})\xrightarrow[n\to\infty]{\text{\tiny dist.}}\ Z$$

\subsection*{Proof of (I) ({\rm as on p. \pageref{I}}) given the {\tt Op-LLT}}
\

By the {\tt LLT} for $(\frak r,\varphi)$ in case ({\tt CN)} and ({\tt Op-LLT}) in case ({\tt NNS)},

\begin{align*}
\sqrt n b(n)^{\kappa}\T^n(1_{A\cap [\varphi_n\in z_n+U,\ \frak r_n\in I+n{E(\frak r)}+r_n]})\ 
\underset{n\to\infty}{\text{\Large $\approx$}}\ 
f_Z(\zeta_n,\rho_n)m(A)|I|m_\Bbb G(U).
\end{align*} for $A\in\mathcal C_\a$  and $I\subset\Bbb R$ an interval so that $A\x\ I\subset Y,\ \&\ U\subset\Bbb G$ a 
compact neighborhood of $0$ with $m_\Bbb G(\bdy U)=0$, 
where $(z_n,r_n)\in\Bbb G\x\Bbb R$ with $z_n=b(n)\cdot\zeta_n,\ \ r_n=\sqrt n\cdot\rho_n$ with $(\zeta_n,\rho_n)=O(1)$.

Fix $t,\ M>0$,\ \ 
then for $n\ge 1$,

\begin{align*}&t=E(\frak r)n+x\sqrt{n}\ \ \text{with} \ \ x=x_{n,t},\ |x|\le ME(\frak r)^{\frac32}\ \ \iff\\ &
  n=\l t -\frac{x'}{{E(\frak r)}^{\frac32}} \sqrt t\ \text{with}\ |x-x'|=O(\tfrac1{\sqrt t})\  \&\ \text{in this case}\\ &
 \T^n(1_{A\cap [\varphi_n\in x(t)+U,\ \frak r_n\in I+t-y]})  \ \sim\ \frac{m(A)|I|}{b(n)^{\kappa}\sqrt n}f_Z(z,x_{n,t})\\ &
\text{as}\ t,\ n\to\infty,\ |x_{n,t}|\le M.
\end{align*}

\

It follows that for fixed $M>0$,
\begin{align*}b(\l t)^{\kappa}I(M,t)&=
b(\l t)^{\kappa} \sum_{n=\l t \pm M\sqrt t}\T^n(1_{A\cap [\varphi_n\in x(t)+U,\ \frak r_n\in I+t-y]})\\ &
\ \underset{t\to\infty}{\text{\Large $\sim$}}\ 
\sum_{n=\l t \pm M\sqrt t}b(n)^{\kappa}\T^n(1_{A\cap [\varphi_n\in x(t)+U,\ \frak r_n\in I+t-y]})
\\ & \ \underset{t\to\infty}{\text{\Large $\sim$}}\   
m(A)|I|m_\Bbb G(U) \sum_{n=\l t \pm M\sqrt t}\frac{f_Z(z,x_{n,t})}{\sqrt{n}}
\end{align*}

Now,
$$\frac1{\sqrt{n}}-\frac1{\sqrt{n+1}}\sim\frac1{2n\sqrt{n}}$$ so
\begin{align*}x_{n,t}-x_{n+1,t}&=\frac{t-{E(\frak r)} n}{\sqrt{n}}-\frac{t-{E(\frak r)} (n+1)}{\sqrt{n+1}}
\\ &=\frac{{E(\frak r)}}{\sqrt{n+1}}+(t-{E(\frak r)} n)(\frac1{\sqrt{n}}-\frac1{\sqrt{n+1}})\\ &
=\frac{{E(\frak r)}}{\sqrt{n}}(1+ O(\frac1{\sqrt n})).
\end{align*}
Thus 

\begin{align*}b(\l t)^{\kappa}\sum_{n=\l t \pm M\sqrt t}&\T^n(1_{A\cap [\varphi_n\in x(t)+U,\ \frak r_n\in I+t-y]})\\ &
\ \underset{t\to\infty}{\text{\Large $\approx$}}\  m(A) |I|m_\Bbb G(U)
\sum_{n=\l t \pm M\sqrt t}(x_{n+1,t}-x_{n,t}) f_Z(z,x_{n,t})\\ &
\xrightarrow[t\to\infty]{}\  m(A) |I|m_\Bbb G(U)\int_{[-M,M]}f_Z(z,x)dx 
\\ &\xrightarrow[M\to\infty]{}  
m(A) |I|m_\Bbb G(U)f_{\mathcal S}(z).
\ \ \CheckedBox\ \text{\tt (I)}\end{align*}

\subsection*{Proof of (II) ({\rm as on p. \pageref{II}}) $\&$ the deviation lemma}
\

 Fix $x:\Bbb R_+\to\Bbb G$ so that $\frac{x(t)}{b(\l t)}\xrightarrow[t\to\infty]{}\ z\in\Bbb R^{\kappa}$ and
 $I,\ U$ as on p. \pageref{flowdisp}, and a compact subset $\frak K:=X\x \mathcal K\x J\subset Z$, then for $(x,z,y)\in\frak K$,
 \begin{align*}II(M,t,(x,z,y))&=
\sum_{n\ge 1,\ |n-\l t|> M\sqrt t}\T^n(1_{A\cap [\varphi_n\in  x(t)-z+V,\ \frak r_n\in I+t-y]})(x)\\ &=
\sum_{n\in A(M,t)}\T^n(1_{A\cap [\varphi_n\in  x(t)-z+V,\ \frak r_n\in I+t-y]})(x)
 \end{align*}
 where 
 $$A(M,t):=\{n\ge 1:\ m([\frak r_n\in t+I-J])>0\ \&\ |n-\l t|> M\sqrt t\}.$$
 Since $\frak r$ is uniformly bounded away from $\{0,\infty\}$, $\exists\ K=K(I,J)>1$ so that
\begin{align*}\tag{\dsrailways}\label{dsrailways}m([\frak r_n\in t+I-J])>0\ \Lra\ n\in [t/K,tK].\end{align*}

 In particular
  \begin{align*}\tag{{\scriptsize\faPaperclip}}\label{Paperclip}
    A(M,t)&\subset B_K(M,t)\\ &:=A(\tfrac{M}K,t)\cap [\tfrac{t}K,tK]\\ &=
    \{n\in [\tfrac{t}K,tK]:\ |n-\l t|> \tfrac{M\sqrt n}K\}.
  \end{align*}

 Thus
   \begin{align*}\tag{\dsjuridical}\label{dsjuridical}
     II(M,t,(x,z,y))\le\sum_{n\in B_K(M,t)}\T^n(1_{A\cap [\varphi_n\in  x(t)-z+V,\ \frak r_n\in I+t-y]})(x).
   \end{align*}

 We'll need the following lemma which follows from a multidimensional version of proposition 4
 in \cite{Pene09} (``extended {\tt LLT} estimate''). Both the lemma and the estimate are proved in the appendix (\S6).
 \proclaim{Deviation Lemma}\label{deviation}
       \

 {\rm (i)}\ There are constants $\G,\ \eta>0$ 
 so that $\forall\ M>0\ \exists\ t_0(M)>0$ such that
  \begin{align*}\T^n1_{A\cap [\frak r_n\in I+t-y]}&\le 
   \G\(\frac{\exp[-\frac{\eta |n-\l t|^2}t]}{\sqrt t}
   +\frac1{n^{\frac32}}\)\\ &
   \forall\ M>0,\ t>t_0(M)\ \&\ n\in B_K(M,t).\end{align*}
    \
    
{\rm (ii)}\ If in addition, $E(\|\v\|^4)<\infty$, then,  possibly changing
$\G,\ \eta,\ t_0(M)$,
we have for $M>0,\ t>t_0(M),\ x(t)\in\Bbb R^\kappa,\ 
\|x(t)\|\le \tfrac{M\sqrt t}{2\l K^2}$:

\begin{align*}\T^n&1_{A\cap [\varphi_n\in  x(t)-z+V,\ \frak r_n\in I+t-y]}\le
   \frac{\G}{t^{\frac{\kappa}2}}\cdot\(\frac{\exp[-\frac{\eta|n-\l t|^2}t]}{\sqrt t}+\frac1{n^{\frac32}}\)\\ &
   \forall\  n\in B_K(M,t).\end{align*}\endproclaim
In the deviation lemma, $(y,z)\in\Bbb R_+\x\Bbb G\ \& I\x V\subset\Bbb R_+\x\Bbb G$ are considered fixed, but  may vary inside a compact subset
of $\Bbb R_+\x\Bbb G$.
\

   \subsection*{Proof of (II) {\rm (as on p. \pageref{II})}  given the Deviation Lemma}
\

   Evidently,
    \begin{align*}\tag{\Wheelchair}\label{Wheelchair}\sum_{n\in B_K(M,t)}\tfrac{1}{n^{\frac32}}\le
    \sum_{n\ge \frac{t}K}\tfrac{1}{n^{\frac32}}\le \tfrac{4\sqrt K}{\sqrt t}
   \end{align*}
 Next we claim that $\exists\ C>0,\ \zeta>0$ so that $\forall\ R,\ t>0$,
 \begin{align*}\tag{\dsbiological}\label{dsbiological}
 \tfrac1{\sqrt t}\sum_{n\in B_K(R,t)}{\exp[-\tfrac{\eta|n-\l t|^2}t]}&\le
 Ce^{-\zeta R^2}.
\end{align*}
\demo{Proof of (\dsbiological)}\ \ Let $N_t:=[\l t]\ \&\ r_t:=\{\l t\}$, then

  \begin{align*}\sum_{n\in B_K(R,t)}\exp[-\tfrac{\eta |n-\l t|^2}t]
&\le
\sum_{n\in B_K(R,t)}\exp[-\tfrac{\eta |(n-N_t)-r_t|^2}t]\\ &\ll                                                                                                                
\sum_{k\in\Bbb Z,\ |k|\ge  \frac{R\sqrt t}K}\exp[-\tfrac{\eta |k-r_t|^2}t]\\ &\ll
2\sum_{k\ge  \frac{R\sqrt t}K}\exp[-\tfrac{\eta k^2}{t}]\\ &\ll 
2\int_{\frac{R\sqrt t}K}^\infty\exp[-\tfrac{\eta y^2}{t}]dy\\ &=
2\sqrt{t}\int_{\frac{R}K}^\infty\exp[-\eta x^2]dy\\ &
\asymp\sqrt t\exp[-\tfrac{\eta R^2}{K^2}]\\ &\ll \sqrt t\exp[-\zeta R^2]
\end{align*}
with $\zeta:=\tfrac{\eta}{2K^2}$.\ \ \Checkedbox\ (\dsbiological).
\demo{Proof in case (i)}
For $(x,z,y)\in\frak K:=X\x K\x J\subset Z$ (compact), using (ii) in the deviation lemma,
 \begin{align*}t^{\frac{\kappa}2}&II(M,t,(x,z,y))\\ &=t^{\frac{\kappa}2}\widehat{\Phi}_t(1_{A\x (x(t)+U)\x I]\cap[|\mathcal N(t)-\l t|\ge M\sqrt t]})(x,z,y)
\\ &= t^{\frac{\kappa}2}\sum_{n\ge 1,\ |n-\l t|> M\sqrt t}\T^n(1_{A\cap [-\varphi_n\in  x(t)-z+U ,\ \frak r_n\in I+t-y]})(x)\\ &\le
t^{\frac{\kappa}2}\sum_{n\in B_K(M,t)}\T^n(1_{A\cap [-\varphi_n\in  x(t)-z+U ,\ \frak r_n\in I+t-y]})(x)\\ &\le
\tfrac{\G}{\sqrt t}\sum_{n\in B_K(M,t)}\exp[-\tfrac{\eta|n-\l t|^2}t]+
\G\sum_{n\in B_K(M,t)}\tfrac{1}{n^{\frac32}}.
 \end{align*}  
 
 Using (\Wheelchair) and (\dsbiological), 
 we obtain
 $$t^{\frac{\kappa}2}II(M,t,(x,z,y))\ll e^{-\zeta M^2}+\tfrac1{\sqrt t}$$
  which suffices for (II).\ \Checkedbox
  
\demo{Proof of (II)  in case (ii)}\ \ 
Define $\frak y:Y\to\Bbb R_+,\ \frak y(x,y):=y$
then using (\dsjuridical) on p. \pageref{dsjuridical} and ({\scriptsize\faAmbulance}) on p. \pageref{faAmbulance} 
 \begin{align*}E(II(M,t,(x,z,\cdot)\|\frak y))&\le\sum_{n\in B_K(M,t)}m(A\cap [\varphi_n\in  x(t)-z+V,\ \frak r_n\in I+t-y])\\ &\le
  \frak c\sum_{n\in B_K(M,t)}m([\varphi_n\in  x(t)-z+V])m([\frak r_n\in I+t-y])
 \end{align*}

with $\frak c$ as in the  weak independence assumption and where $E(\cdot\|\frak y)$ denotes conditional expectation on $(Y,\nu)$ 
with respect to $\frak y$.

As $t\to\infty$ and $n=n(t)\in B_K(M,t)$ (as on p. \pageref{Paperclip}), $n(t)=K^{\pm 1}t$ and $b(\l n(t))\asymp b(t)$.
\

Thus by the {\tt LLT} for $(X,m,T,\v),\ \exists\ \frak k>1$ so that
$$b(t)^\kappa m([\varphi_n\in  x(t)-z+V])\le \frak k\ \forall\ t>0\ \text{large},\ n\in B_K(M,t).$$

By part (i) of the deviation lemma,
$$m([\frak r_n\in I+t-y])\le 
\G\(\frac{\exp[-\frac{\g|n-\l t|^2}t]}{\sqrt t}+\frac1{n^{\frac32}}\)\ \forall\ t>0,\ n\in B_K(M,t).$$
It follows that

\begin{align*}
 b(t)^\kappa E(II(M,t,&(x,z,\cdot)\|\frak y))\le 
 \frak c \G \frak k\sum_{n\in B_K(M,t)}(\tfrac{
 \exp[-\frac{\g|n-\l t|^2}t]}{\sqrt n}+\tfrac1{n^{\frac32}})\\ &
 \ll e^{-\zeta M^2}+\tfrac1{\sqrt t}\end{align*}
as before by (\dsbiological) $\&$ (\Wheelchair). This suffices for (II).\ \Checkedbox

\subsection*{ Proof of theorem 2}
\

Under the assumptions for theorem 2, we have by \cite{AD} that the measure preserving transformation
$(X\x\Bbb G,m\x m_\Bbb G,T_\varphi)$ is pointwise dual ergodic with 
$$a_n(T_\v)\sim \sum_{k=1}^n\frac{f_\mathcal S(0)}{b(k)^{\kappa}}.$$ 
\

Since $(X\x\Bbb G,m\x m_\Bbb G,T_\varphi)$ is a {\tt good section} for the  skew product semiflow $(Z,\mu,\Phi^{(J)})$ in the sense of \cite{mode},
we have (as in \cite{nice})  that each $\Phi^{(J)}_t$ is pointwise dual ergodic with 
$$a_n(\Phi^{(J)}_t)\sim \sum_{k=1}^n\frac{f_\mathcal S(0)}{b(\l tk)^{\kappa}}.$$

Now let 
 $A\in\mathcal C_\a$,
$I\subset\Bbb R$ be an interval and $U\subset\Bbb G$ be a compact neighborhood with $m_{\Bbb G}(\bdy U)=0$.
\

By (I), as $t\to\infty$,
\begin{align*}\widehat{\Phi^{(J)}}_t(&1_{A\x U\x I})(\om,z,y)=\\ &=
\widehat{\Phi}_t(1_{(A\x I)\cap [J^{(\varphi)}(t,\cdot)\in  z-U]})(\om,y)\\ &
\gtrsim\ f_{\mathcal S}(0)\mu(A)\text{\tt Leb}\,(I)m_\Bbb G(U)\cdot\frac1{b(\l t)^{\kappa}}.
\end{align*}
By pointwise dual ergodicity, for fixed $t>0$, as $N\to\infty$,
\begin{align*}\sum_{k=1}^N\widehat{\Phi}_{kt}(&1_{(A\x I)\cap [J^{(\varphi)}(t,\cdot)\in  U]})(\om,y)\\ &
\sim\ f_{\mathcal S}(z)m(A)\text{\tt Leb}\,(I)m_\Bbb G(U)\sum_{k=1}^N\frac1{b(\l kt)^{\kappa}}.
\end{align*}
The rest of the proof
  follows the proof of proposition 4.2 in \cite{nice}.\ \Checkedbox
   \
   
   The proof of theorem 2 actually establishes the (possibly) more general 
   \proclaim{Abstract Proposition}\label{abstract}
 \
 
 If  $(X,m,T,\alpha)$ is probability preserving,  mixing,
 {\tt G-M}-- or   {\tt AFU}-- map and 
  $(\varphi,\frak r):X\to\Bbb G\x\Bbb R_+$, with $E(\frak r)<\infty$, 
 satisfies   $\mathcal D^{(T)}_\alpha((\varphi,\frak r))<\infty$ and the operator local limit theorem 
 corresponding to the distributional convergence 
 $$(\tfrac{\frak r_n-n{E(\frak r)}}{\b(n)},\tfrac{\varphi_n}{b(n)})\xrightarrow[n\to\infty]{\text{\tiny dist.}}\ Z$$
 (as on  p. \pageref{Op-LLT}) 
 where $b(n)$ is $p$-regularly varying with $0<p\le 2$, $\b(n)$ is $q$-regularly varying with $1<p\le 2$
 and $Z$ is a globally supported $\rv$ on $\Bbb R^{\kappa+1}$ admitting a  smooth probability density function; then the skew product semiflow is conditionally, rationally weakly mixing.\endproclaim
 
 \
 
 \subsection*{Remark}\ It seems reasonable to expect the rational weak mixing of:
 \sbul the geodesic flow on the unit  tangent bundle of $\Bbb C\setminus\Bbb Z$ (as in 
 \cite{mode}): and 
 \sbul the planar Lorentz systems considered in \cite{Sz-V1} $\&$ \cite{Sz-V2};
 
 \
 
as these flows are (measure theoretic) natural extensions of skew product semiflows for which the above abstract proposition might hold.
  \subsection*{ Proof of theorem 3}
   \par Let 
$V$ be the $\Bbb Z$-cover of the compact hyperbolic surface $M$ (as on p. \pageref{geodesic}) 
 equipped with a covering map $p:V\to M$ so that
there exists a monomorphism $\g:\Bbb Z\to\text{\tt Isom}\,(V)$, such that
for $y\in V,\ p^{-1}\{p(y)\}=\{\g(n)y:\ n\in\Bbb Z\}$. 
\

The corresponding tangent map, also denoted  $p:U(V)\to U(M)$ is equivariant with the geodesic flows. 

\par As $g^M$ is an Anosov flow, by Bowen's theorem (\cite{Bowen}, see also \cite{PaPo}), there is a suspended flow
$(X,m,T)^\frak r$ (see p. \pageref{faLinux}) and  a continuous, measure theoretic isomorphism 
$\pi:(X,m,T)^\frak r\to (U(M),\L,g)$.
\

Here, $X\subset S^\Bbb Z$ with $S$ finite, is a mixing, two-sided SFT with $T$ the shift and  $m$  a Gibbs measure.
\

Denote the corresponding one-sided SFT
by $(X_+,m_+,T_+)$ with
 $X_+:=\{\varpi(x):\ x\in X\}$ where  $\varpi:X\to S^\Bbb N,\ \varpi(x):=x|_\Bbb N;$  $m_+:=m\circ\varpi^{-1}\ \&\ T_+$ is the shift.
\

The roof function $\frak r:X\to\Bbb R_+$ is bounded away from $\{0,\infty\}$ and {\it one-sided} in the sense that
$\frak r=\frak r_+\circ\varpi$ where  $\log\frak r_+:X_+\to\Bbb R$ is H\"older continuous:
$$|\log\frak r_+(x)-\log\frak r_+(y)|\le M\th^{t(x,y)}$$
where $M>1,\ \th\in (0,1)\ \&\ t(x,y):=\min\,\{n\ge 1:\ \ x_n\ne y_n\}\le\infty$.
\

As in \cite{Rees}, $\exists$ a block map $\phi:X\to\Bbb Z$ and a continuous,  measure theoretic isomorphism 
$$\frak P:(X\x\Bbb Z,m\x\#,T_\phi)^{\overline{\frak r}}\to (U(V),\L,g)$$ so that
$\frak P\circ Q_n=\g(n)\circ\frak P$. 
Here $(X\x\Bbb Z,m\x\#,T_\phi)^{\overline{\frak r}}$ is as in (\faLinux) on p. \pageref{faLinux}.
\

Recall that, as mentioned on p. \pageref{topflow}, $(X,m,T)^\frak r=(Y,\nu,\psi)$ is a continuous flow on a compact metric space and 
$(X\x\Bbb Z,m\x\#,T_\phi)^{\overline{\frak r}}$  is  a continuous flow on the 
locally compact, Polish, metric space. 	The continuity of $\frak P$ above is with respect to this topology.
\

It suffices to prove ({\scriptsize\faBus}), ({\scriptsize\faShip}) $\&$ (\Yinyang) for $(X\x\Bbb Z,m\x\#,T_\phi)^{\overline{\frak r}}$.
\

To this end, we show first 
\Par \ \ $\phi:X\to\Bbb Z$ is centered and  $(\frak r,\phi):X\to\Bbb R\x\Bbb Z$ is aperiodic.
\pf \ 
By \cite{Rees}, $(U(V),\L,g)$ is conservative, whence (\cite{Hopfbook}), ergodic. These properties pass to
$(X\x\Bbb Z,m\x\#,T_\phi)^{\overline{\frak r}}$, whence $E(\phi)=0$.

\ As in  \cite{Sharp}, $V$ is {\it homologically full} in the sense that as $t\to\infty,\ \exists$ 
exponentially many closed geodesics of length $\le t$ in each homology class.
\

\sms Therefore,  by the lemma in \cite{Solomyak}  the function $(\frak r,\phi):X\to\Bbb R\x\Bbb Z$ is   aperiodic.\ \ \Checkedbox\ \P

 \

Next, by Sinai's theorem (\cite{Sinai}, see also \cite{Sarig}) there are  block functions
$\xi:X\to\Bbb Z\ \&\ ,\ \phi^+:X_+\to\Bbb Z$ so that
$\phi=\xi-\xi\circ T+\phi^+.$
\

It follows that
$$(X\x\Bbb Z,m\x\#,T_\phi)^{\overline{\frak r}}\ \ \&\ \ (X\x\Bbb Z,m\x\#,T_{\phi^+\circ\varpi})^{\overline{\frak r_+\circ\varpi}}$$
are {\tt RI} conjugate.
\

Thus,  it 
 suffices to prove ({\scriptsize\faBus}), ({\scriptsize\faShip}) $\&$ (\Yinyang) for 
 $$(\widetilde{Z},\widetilde{\mu},\Psi):=(X\x\Bbb Z,m\x\#,T_{\phi^+\circ\varpi})^{\overline{\frak r_+\circ\varpi}}.$$

 \

The continuous natural extension  $\varpi:(X,m,T)\to (X_+,m_+,T_+)$ induces a continuous natural extension 
$$\Pi:(\widetilde{Z},\widetilde{\mu},\Psi)\to 
(Z,\mu,\Phi):=(X_+\x\Bbb Z,m_+\x\#,T_{\phi^+})^{\overline{\frak r_+}}$$

Let
$$\mathcal H:=L_T\otimes C_c(\Bbb R_+)\otimes C_c(\Bbb Z)\subset C_c(Z).$$
By theorem 1(i), with  $b(t)\propto\sqrt t$ 
and some  centered Gaussian random variable  $\mathcal S$  on $\Bbb R$, 
for each $A\in\mathcal H$ and 
$$k:\Bbb R_+\to\Bbb G,\ \tfrac{k(t)}{b(t)}\xrightarrow[t\to\infty]{} z,\ \sup_{t>0}\tfrac{\|k(t)\|}{b(t)}<\infty,$$
\begin{align*}\tag{{\scriptsize\faRocket}}b(t)\widehat{\Phi}_t(A\circ Q_{k(t)})(x)&\xrightarrow[t\to\infty]{}f_\mathcal S(z)\mu(A)
\\ & \text{uniformly on compact subsets of}\ Z.\end{align*}
Applying this with $k\equiv 0$ repeatedly, it follows that for $J\ge 1,\ A_i\in\mathcal H$,  ($1\le i\le J$), as $t_1,\dots,t_J\to\infty$:
\begin{align*}\tag{{\scriptsize\faBomb}}\prod_{j=1}^Jb(t_j)&\widehat{\Phi}_{t_J}(A_J\widehat{\Phi}_{t_{J-1}}(A_{J-1}\cdots
\widehat{\Phi}_{t_{2}}(A_{2}\widehat{\Phi}_{t_{1}}(A_{1})\cdots)\\ &
\rightarrow{}\prod_{j=1}^Jf_\mathcal S(z_j)\mu(A_j)\ \  \text{uniformly on compact subsets of}\ Z.\end{align*}
 From   ({\scriptsize\faRocket}), we calculate that
 \begin{align*}\tag{{\scriptsize\faFlash}}b(t)\widetilde{\mu}(U\cdot V\circ &Q_{k(t)}\circ\Psi_t)\xrightarrow[t\to\infty,\ \frac{k(t)}{b(t)}\to z]{}f_\mathcal S(z)
 \widetilde{\mu}(U)\widetilde{\mu}(V) \\ &\forall\  U,\ V\in\widetilde{\mathcal H}:=\{G\circ\pi\circ\Psi_{-s}:\ G\in\mathcal H,\ s>0\}
 \end{align*}
 and from ({\scriptsize\faBomb}) we see that for $N\ge 1,\  A,\ A_1,\dots,A_N\in\widetilde{\mathcal H}$,  
\begin{align*}\tag{{\scriptsize\faFire}}
\L(A\cap\prod_{i=1}^{N}A_i\circ\Psi_{t_i})\ \sim\ \L(A)\prod_{i=1}^{N}\tfrac{\L(A_i)}{b(\D t_i)}\end{align*}
as $\u t=\{0=t_0<t_1<\dots<t_N\},\ \D t_i:=t_i-t_{i-1}\to\infty$, ($1\le i\le N$).
\

Evidently, $\widetilde{\mathcal H}\subset C_c(\widetilde{Z})$ is closed under multiplication. Moreover $\widetilde{\mathcal H}$ separates points of $\widetilde{Z}$.
Thus by the Stone-Weierstrass theorem, $\text{\tt Span}\,\widetilde{\mathcal H}$ is locally dense in $C_c(\widetilde{Z})$ and
 $\widetilde{\mathcal H}\subset C_c(\widetilde{Z})$ is indeed a measure determining collection on $\widetilde{Z}$.
\

Standard, repeated application of the weak convergence lemma now establishes
\bul  ({\scriptsize\faFlash}) $\forall\  U,\ V\in C_c(\widetilde{Z})$ and
\bul ({\scriptsize\faFire}) $\forall\ N\ge 1,\  A,\ A_1,\dots,A_N\in C_c(\widetilde{Z}).$
\

As remarked above, this proves theorem 3(o) $\&$ (i).\ \ \Checkedbox
\

To conclude, it follows from ({\scriptsize\faRocket}) with $k(t)=0$ and proposition 3.1 in \cite{nice} that
$(Z,\mu,\Phi^{(J)})$ is rationally, weakly mixing of order $2$ as in (\Yinyang) (on p. \pageref{Yinyang}).
\

This property passes to its measure theoretic natural extension $(U(V),\L,g)$.
this proves theorem 3(ii).\ \ \Checkedbox

\section*{\S4 Operator Local limit theorem}\ \ 
 \

 We'll use the spectral and perturbation  theories of the transfer operator of a fibered system $(X,m,T,\alpha)$ which is assumed to be
 mixing, probability preserving and either a {\tt G-M} map, or an {\tt AFU} map. Details can be found in \S2 in \cite{AD} for  
 {\tt G-M} maps and \S5 in \cite{ADSZ} for   {\tt AFU} maps.
 \
 
 Let $\Bbb G\le\Bbb R^{\kappa},\ $, a subgroup of full dimension and let $\varphi=(\varphi^{(1)},\dots,\varphi^{(\kappa)}):\ X\to\Bbb G$ satisfy
$$\mathcal D^{(T)}_\alpha(\varphi)<\infty$$
where $\mathcal D^{(T)}_\alpha$ is as on p. \pageref{osc}.
\

For $t\in\widehat{\Bbb G}$ define $P_t:L^1(m)\to L^1(m)$ by
$$P_t(f):=\T(e^{i\langle t,\varphi\rangle}f),$$

then, each $P_s:L_T\to L_T$ acts quasicompactly where as above $L_T:=H_\th$  in the {\tt G-M} case $\&$ $L_T=\text{\tt BV}$ in the {\tt AFU} case.

\par Moreover $s\mapsto P_s$ is continuous ($\widehat{\Bbb G}\to \text{\tt Hom}(L_T,L_T)$).

For small $t\in\widehat{\Bbb G}$, the characteristic function operator $P_t$ has a simple, dominant eigenvalue and indeed,
by Nagaev's theorem:
\sms 1) There are constants $\e>0,\ K>0$ and $\th\in (0,1)$; and  continuous functions
$\l:B(0,\e)\to B_{\Bbb C}(0,1),\  N:B(0,\e)\to \text{\tt Hom}(L,L)$ such that
$$\|P_t^nh-\l(t)^n N(t)h\|_L\le K\th^n\|h\|_L
 \ \ \ \ \forall\ |t|<\e,\ n\ge 1,\ h\in L$$
where $\forall |t|<\epsilon$, $N(t)$ is a projection onto a one-dimensional subspace (spanned by $g(t):=N(t)\mathbb{1}$).
\sms  2) If $\varphi$ is aperiodic, then
$\forall\ \widetilde M>0,\ \e>0,\
\exists\ K'>0$ and $\th'\in (0,1)$ such that
$$\|P_t^nh\|_L\le K'\th^{\prime n}\|h\|_L\
\ \ \forall\ \e\le\ |t|\le \widetilde M,\ h\in L.$$

 \subsection*{Operator Local limits}

Now suppose that $0<p<2$ and that the  distribution  of $\varphi$ is
in the strict domain of attraction of
a symmetric $p$-stable random variable $\mathcal S$; equivalently
for some $\frac1p$-regularly varying sequence $(b(n):\ n\ge 1)$,
$$E(\exp[i\<\tfrac{\varphi}{b(n)},u\>])^n\xrightarrow[n\to\infty]{}\ E(e^{i\<\mathcal S,u\>}).$$
As shown in \cite{AD}, and \cite{ADSZ}, 
\begin{align*}\tag{{\scriptsize\faPlug}}\label{faPlug}
 E(\exp [i\<\tfrac{\varphi_n}{b(n)},u\>])^n\approx \l(\tfrac{u}{b(n)})^n
 \xrightarrow[n\to\infty]{}\ E(e^{i\<\mathcal S,u\>})
\end{align*}
Using this, one obtains  (by the methods of
\cite{Breiman}\ $\&$\ \cite{Stone}) the {\tt conditional LLT}
\begin{align*}
b(n)^{\kappa}\T^n(1_{A\cap [\varphi_n\in z_n+U]})\ \underset{n\to\infty}{\text{\Large $\approx$}}\ 
f_\mathcal{S}(\zeta_n)\mu(A)m_\Bbb G(U).
\end{align*} for $A\in\mathcal C_\a$   $\&\ U\subset\Bbb G$ a 
compact neighborhood of $0$ with $m_\Bbb G(\bdy U)=0$, 
where $z_n\in\Bbb G$ with $z_n=b(n)\cdot\zeta_n,\ \zeta_n=O(1)$.

\demo{Proof of the {\tt Op-LLT} {\small\rm as on p. \pageref{Op-LLT}  given the {\tt Infinite Divisibility Lemma}}}

We have  
\begin{align*}\tfrac{\frak r_n-n{E(\frak r)}}{\sqrt n}\xrightarrow[n\to\infty]{\text{\tiny dist.}}\mathcal G\ \&\ \tfrac{\varphi_n}{b(n)}\xrightarrow[n\to\infty]{\text{\tiny dist.}}\mathcal S\end{align*}

 where $\mathcal G$ is centered Gaussian on $\Bbb R$ and $\mathcal S$ is symmetric, $p$-stable on $\Bbb R^{\kappa}$
 \
 
 It follows that $((\tfrac{\frak r_n-n{E(\frak r)}}{\sqrt n},\tfrac{\varphi_n}{b(n)}):\ n\ge 1)$ is a tight sequence of random variables on $\Bbb R\x\Bbb R^{\kappa}$.
 \
 
 Let 
$Z=(\mathcal G,\mathcal S)\in\text{\tt RV}\,(\Bbb R\x\Bbb R^{\kappa})$ be a weak limit point of the sequence. 
 By the {\tt Infinite Divisibility Lemma} (see p. \pageref{inf-div}) $Z$ is infinitely divisible. By the Levy-Ito decomposition (see \cite{Sato}),
 $Z=G+C+c$ where $G,\ C\in\text{\tt RV}\,(\Bbb R\x\Bbb R^{\kappa})$ with $G$  centered  Gaussian, $C$  compound Poisson 
 (e.g. $p$-stable with $p<2$) and $c$ is constant. This decomposition is unique and, moreover, $G\ \&\ C$  are independent.
 \
 
Now, $Z=(\mathcal G,0)+(0,\mathcal S)$, whence the Levy-Ito decomposition  is   $G=(\mathcal G,0)$,  $C=(0,\mathcal S)$ and $c=0$, with the 
conclusion that $\mathcal G\ \&\ \mathcal S$ are independent.
\

 This determines $Z$ uniquely and so
 $$(\tfrac{\frak r_n-n{E(\frak r)}}{\sqrt n},\tfrac{\varphi_n}{b(n)})\xrightarrow[n\to\infty]{}\ Z.$$
It follows from the analog of ({\scriptsize\faPlug}) (on p. \pageref{faPlug})  that if $\l(x,y)$ is the dominant eigenvalue of 
 $$P_{(x,y)}f=\T(\exp[ix(\frak r-{E(\frak r)})+i\<y,\varphi\>]f),$$ then
 $$\l(\tfrac{x}{\sqrt{n}},\tfrac{y}{b(n)})^n\xrightarrow[n\to\infty]{}\ E(e^{ix\mathcal G})E(e^{i\<y,\mathcal S\>})$$
 uniformly on compact subsets.
 \

 Now, the stable characteristic functions have the forms
 $$E(e^{ix\mathcal G})=e^{-ax^2}\ \&\ 
E(e^{i\langle y,\mathcal S\rangle})=e^{-c_{p,\nu}(y)}.$$
Here $a>0\ \&\ c_{p,\nu}(y):=\int_{S^{\kappa-1}}|\<y,s\>|^p\nu(ds)$
where $\nu$ is a symmetric
measure on $S^{\kappa-1}$.
 In other words,
  \begin{align*}\tag{\Bicycle}\label{Bicycle}
    -n\log\l(\tfrac{x}{\sqrt{n}},\tfrac{y}{b(n)})\xrightarrow[n\to\infty]{}\ ax^2+c_{p,\nu}(y) 
  \end{align*}

     uniformly on compact subsets.
\

 As in \cite{Stone}, fix $f\in L^1(\Bbb R\x\Bbb G)$ so that $\widehat{f}:\widehat{\Bbb R\x\Bbb G}\to\Bbb C$  is continuous, 
 compactly supported,  $A\in\mathcal C_\a$ and $(r_n,z_n)\in\Bbb R\x\Bbb G$ with 
 $z_n=b(n)\cdot\zeta_n,\ \ r_n=\sqrt n\cdot\rho_n$ with $(\zeta_n,\rho_n)\to (\zeta,\rho)$.
 
 By Nagaev's theorem, (\Bicycle) (p. \pageref{Bicycle})  and ({\tt J-Ap}) as on p. \pageref{J-Ap}   (respectively), there exist\ $\vartheta\in (0,1)\ \&\ \d>0$ so that 
   {\small\begin{align*}&\tag{i}\sup_{(x,y)\in B(0,\d)}\,\|P_{(x,y)}^n1_A-\l(x,y)^nN(x,y)1_A\|_{\text{\tt\tiny Hom}\,(L,L)}=O(\vartheta^n)\ \text{as}\ n\to\infty,
    \\ &\tag{ii}  -n\text{\tt Re}\log\l(\tfrac{x}{\sqrt{n}},\tfrac{y}{b(n)})\ge \th(ax^2+c_{p,\nu}(y))\ \forall\ n\ge 1,\ (x,y)\in B(0,\d),
    \\ &\tag{iii} \sup_{(x,y)\in \text{\tt\tiny supp}\,\widehat{f}\setminus B(0,\d)}\,\|P_{(x,y)}^n1_A\|_L=O(\vartheta^n)\ \ \text{as}\ n\to\infty.
   \end{align*}}
 \

 We have, using (i) and (iii):
 \begin{align*}
  \sqrt n& b(n)^{\kappa}\T^n(1_Af(\frak r_n-{E(\frak r)} n-r_n,\varphi_n-z_n))\\ & =
  \sqrt n b(n)^{\kappa}\int_{\widehat{\Bbb R\x\Bbb G}}e^{-ir_nx-i\<y,z_n\>}\widehat{f}(x,y)P_{(x,y)}^n1_Adxdy\\ &=
   \sqrt n b(n)^{\kappa}(\int_{B(0,\d)}+\int_{\widehat{\Bbb R\x\Bbb G}\setminus B(0,\d)})e^{-ir_nx-i\<y,z_n\>}\widehat{f}(x,y)P_{(x,y)}^n1_Adxdy\\ &=
  \sqrt n b(n)^{\kappa}\int_{B(0,\d)}e^{-ir_nx-i\<y,z_n\>}\widehat{f}(x,y)P_{(x,y)}^n1_Adxdy+O(\sqrt n b(n)^{\kappa}\vartheta^n)\\ &=
  \sqrt n b(n)^{\kappa}\int_{B(0,\d)}e^{-ir_nx-i\<y,z_n\>}\widehat{f}(x,y)\l(x,y)^nN(x,y)1_Adxdy+O(\sqrt n b(n)^{\kappa}\vartheta^n).
   \end{align*}
Writing 
$\D_n(x,y):=(\tfrac{x}{\sqrt n},\ \tfrac1{b(n)}\cdot y)$ and changing variables, 
 \begin{align*}\sqrt n& b(n)^{\kappa}\int_{B(0,\d)}e^{-ir_nx-i\<y,z_n\>}\widehat{f}(x,y)\l(x,y)^nN(x,y)1_Adxdy\\ &=
 \int_{\Bbb R^{\kappa+1}}1_{\D_n^{-1}B(0,\d)}(x,y)e^{-i\rho_nx-i\<y,\zeta_n\>}\widehat{f}(\D_n(x,y))
 \l(\D_n(x,y))^nN(\D_n(x,y))1_Adxdy\\ &\xrightarrow[n\to\infty]{}\ 
\widehat{f}(0)\cdot N(0)1_A\cdot\int_{\Bbb R\x\Bbb R^{\kappa}}e^{-i\rho x-i\<y,\zeta\>}e^{-ax^2-c_{\nu,p}(y)}dxdy\\ &=
 \int_\Bbb G f(y)dy\cdot m(A)\cdot f_Z(\rho,\zeta).
   \end{align*}
The convergence here is by Lebesgue's dominated  convergence theorem, since by (ii):
$$1_{\D_nB(0,\d)}(x,y)|\l(\D_n(x,y))|^n\le \exp[-\th(ax^2+c_{p,\nu}(y))].$$
({\tt Op-LLT}) follows from this as in \cite{Breiman}.\ \ \Checkedbox

\f{\bf Remark}
\

As in the proof of lemma 6.4 in \cite{AD}, ``changing variables'' in  (\Bicycle), we see that 
\begin{align*}\tag{{\scriptsize\faPuzzlePiece}}\label{Puzzle}-\log\l(su,tv)\sim as^2+\frac{c_{p,\nu}(v)}{b^{-1}(\frac1{t})}\ \ \text{as}\ \ \|(s,t)\|\to 0,\ 
  b(\tfrac1t)\propto\frac1{s^2}  \end{align*}
  uniformly in $u=\pm 1,\ v\in\Bbb S^{\kappa-1}$. {\it c.f.} lemma 2.4 in \cite{Thomine}.

\section*{\S5 Infinite divisibility lemma}
In this section, we complete the proof of the {\tt Op-LLT} by establishing:
\proclaim{Infinite divisibility lemma}\label{inf-div}
\

Let $(X,m,T,\alpha)$ be a reverse $\phi$-mixing (as on p. \pageref{fatso}) 
fibered system and let $\varphi:X\to\Bbb G$ be $\th$-H\"older continuous
on each $a\in\alpha$ with $D_{\alpha,\th}(\varphi)<\infty$. 
\

Suppose that $\D_n=\text{\tt diag}(\d_1^{(n)},\dots,\d_\kappa ^{(n)})$ are $\kappa \x \kappa $ diagonal matrices satisfying
\sms {\rm(i)}\  each $\d_j^{(n)}=\tfrac1{b_j(n)}$ with 
 $b_j(n),\ \g_j$-regularly varying with $\g_j>0$ 

 and 
\sms {\rm (ii)} $\{\D_n\varphi_n:\ n\ge 1\}$ is uniformly tight.
\

Let   $n_k\to\infty$ and suppose that
$$\D_{n_k}\varphi_{n_k}\xrightarrow[k\to\infty]{\text{\tiny dist.}}\ Z$$
where $Z$ is a random variable on $\Bbb R^{\kappa}$, then
$Z$ is infinitely divisible.
\endproclaim
\demo{Proof}
\

By tightness,  $c(t):=\sup_{n\ge 1}m([\|\D_n\varphi_n\|>t])\xrightarrow[t\to\infty]{}\ 0$.

\

Suppose that $1\le r\le n$, then 
$$\|\D_n\varphi_r\|=\|\D_n\D_r^{-1}\D_r\varphi_r\|\le \|\D_n\D_r^{-1}\|\cdot\|\D_r\varphi_r\|$$
whence for $\e>0$,
\begin{align*}\tag{$\diamondsuit$}\label{diamondsuit} m([\|\D_n\varphi_r\|>\e])&\le m([\|\D_r\varphi_r\|>\tfrac{\e}{\|\D_n\D_r^{-1}\|}])\\ &\le
c(\tfrac{\e}{\|\D_n\D_r^{-1}\|}])\\ &\to 0\ \text{as}\ n\ge r\to\infty,\ \|\D_n\D_r^{-1}\|\to 0. 
\end{align*}
Next, for $Z$ a random variable, let
$$\s(Z):=E\left(\frac{\|Z\|}{1+\|Z\|}\right),$$ then for $Z_n\ (n\ge 1)$ random variables
\begin{align*}&\s(\sum_{k=1}^nZ_k)\le \sum_{k=1}^n\s(Z_k),\\ &
Z_n\xrightarrow[n\to\infty]{P} 0\ \ \iff\ \ \s(Z_n)\xrightarrow[n\to\infty]{} 0. 
\end{align*}
Similar to ($\diamondsuit$), we have
\begin{align*}\tag{$\spadesuit$}\label{spadesuit}  C(t):=\sup\,\{\s(\D_n\varphi_r):\ \ n\ge r\ge 1,\ \|\D_n\D_r^{-1}\|\le t\}\xrightarrow[t\to 0+]{}\ 0.
\end{align*}\demo{Proof of ($\spadesuit$)}\ \ If not, then $\exists\ \nu_k\ge r_k\ge 1$ and $\e>0$ so that
$$\|\D_{\nu_k}\D_{r_k}^{-1}\|\to 0\ \&\ \s(\D_{\nu_k}\varphi_{r_k})\ge \e>0.$$
But in this case, by ($\diamondsuit$)\ $\D_{\nu_k}\varphi_{r_k}\xrightarrow[k\to\infty]{P}0$, whence
$\s(\D_{\nu_k}\varphi_{r_k})\xrightarrow[k\to\infty]{} 0$. Contradiction.\ \ \Checkedbox ($\spadesuit$)

\

Using the regular variation of each $b_j$ and possibly passing to a subsequence of  $n_k\to\infty$, we can ensure that 
\begin{align*}\tag{\dsmilitary}\ b_j(n_k)=2^{\pm 1}k^{\g_j}b_j(m_k)\ \forall\ 1\le j\le \kappa\ \text{where}\ m_k:=\lfl\tfrac{n_k}k\rfl\end{align*}
in addition to
\begin{align*}\tag{\dsmathematical}\label{dsmathematical} k\max_{1\le J\le 2k}\|\D_{n_k}\D_{J}^{-1}\|\to 0                                               
                                              \end{align*}

and
$$\D_{n_k}\varphi_{n_k}\xrightarrow[k\to\infty]{\text{\tiny dist.}}\ Z.$$

 Write $n_k=km_k+r_k$ with $0\le r_k<m_k$.
 \

Let $b:=\min_j\g_j>0$. By (\dsmilitary), we have  $\|\D_{n_k}\D_{m_k}^{-1}\|\le \frac2{k^b}\to 0$ and hence by  ($\spadesuit$), that
\begin{align*}\D_{n_k}\varphi_{m_k}\xrightarrow[k\to\infty]{P}\ 0. 
\end{align*}

Next, we claim that
\begin{align*}\tag{\Radioactivity}\label{Radioactivity}    \D_{n_{k}}\sum_{j=1}^{k-1}\varphi_{2k}\circ T^{jm_k-2k}\xrightarrow[k\to\infty]{P}\ 0.
\end{align*}
\demo{Proof of (\Radioactivity)} 

\begin{align*}\s\(\D_{n_{k}}\sum_{j=1}^{k-1}\varphi_{2k}\circ T^{(j+1)m_k-2k}\)&\le k \s(\D_{n_{k}}\varphi_{2k})
\\ &\xrightarrow[k\to\infty]{}0\ \ \text{by (\dsmathematical).\ \ \Checkedbox(\Radioactivity)}
\end{align*}
Next, we'll use (\Radioactivity) to show :

\begin{align*}\tag{\dschemical} \label{dschemical}          \D_{n_k}\sum_{j=0}^{k-1}\(\varphi_{m_k-2k}\circ T^{jm_k}+\varphi_{r_k}\circ T^{jm_k}\)
\xrightarrow[k\to\infty]{\text{\tiny dist.}}Z. 
\end{align*}
\demo{Proof}
\begin{align*}\D_{n_k}\sum_{j=0}^{k-1}\(\varphi_{m_k-2k}\circ T^{jm_k}&+\varphi_{r_k}\circ T^{jm_k}\)-\D_{n_k}\varphi_{n_k}\\ &=
\D_{n_{k}}\sum_{j=1}^{k-1}\varphi_{2k}\circ T^{jm_k-2k}\\ &\xrightarrow[k\to\infty]{P}\ 0\ 
\text{by (\Radioactivity). This proves (\dschemical).\ \Checkedbox}
\end{align*}
\

 For $n\ge 1,\ a\in\alpha_n$, fix $\zeta_n(a)\in T^na$ so that $\exists\ z_n(a)\in a$ satisfying
 $T^nz_n(a)=\zeta_n(a).$
 \
 
Now define $\pi_n:X\to X$ by 
$$\pi_n(x):=z_n(\alpha_n(x)).$$
It follows that for $n,\ k\ge 1$,
\begin{align*}
 \tag{\dstechnical}\label{dstechnical}  |\varphi_n-(\varphi\circ\pi_{n+k})_n|\le\frac{D_{\alpha,\th}(\varphi)\th^k}{1-\th}.
\end{align*}

\

Set 
$$W_{k,j}:=\D_{n_k}\varphi_{m_k-2k}\circ T^{jm_k},\ 0\le j\le k-1\ \&\ W_{k,k}:=\D_{n_k}\varphi_{r_k}\circ T^{km_k}.$$
$$Y_{k,j}:=\D_{n_k}\varphi_{m_k-2k}\circ\pi_{m_k-k}\circ T^{jm_k},\ 0\le j\le k-1\ \&\ Y_{k,k}:=\D_{n_k}\varphi_{r_k}\circ T^{km_k}.$$
By (\dstechnical), for each $0\le j\le k$,
$$\|Y_{k,j}-W_{k,j}\|\le \frac{\|\D_{n_k}\|D_{\alpha,\th}(\varphi)k\th^k}{1-\th}$$ whence
\begin{align*}\sum_{j=0}^k\|Y_{k,j}-W_{k,j}\|\le \frac{\|\D_{n_k}\|D_{\alpha,\th}(\varphi)k\th^k}{1-\th}\xrightarrow[k\to\infty]{}0. 
\end{align*}
It follows that 
\begin{align*}&\tag{i}\sum_{j=0}^kY_{k,j}\xrightarrow[k\to\infty]{\text{\tiny dist.}}Z \ \ \&\\ &
\tag{ii}\max_j\,\s(Y_{k,j})\xrightarrow[k\to\infty]{} 0. 
\end{align*}

Now $\{(Y_{k,k-j}:\ 0\le j\le k):\ k\ge 1)$ is a $\phi$-mixing, triangular array as in \cite{Samur} and as shown there,
$Z$ is infinitely divisible.\ \ \Checkedbox

\

 \section*{\S6 Appendix}
  In this appendix, we state and prove an Extended {\tt LLT} estimate and deduce the 
  the Deviation lemma from it.
  \
  
  The estimate is a multidimensional version of proposition 4 in \cite{Pene09}. For earlier work on the iid case, see
  \cite[Chapter VII,\S3]{Petrov}. 
 
 \
 
  Let 
\bul $(X,m,T,\a)$ be a mixing fibered system, either {\tt GM} or {\tt AFU};
\bul $d\ge 1,\ \Bbb H\le\Bbb R^d$ be a closed subgroup of full dimension; 
\bul  $\Psi:X\to\Bbb H$ be aperiodic, $D_\a^{(T)}(\Psi)<\infty$ and $E(\|\Psi\|^4)<\infty$.

\

\

\proclaim{Extended {\tt LLT} estimate\ \ (\cite{Pene09})}
\

Fix $U\subset\Bbb H,\ m_\Bbb H(\bdy U)=0$ precompact. There exist
 $\Gamma>0,\ \eta>0$ 
 so that
   \begin{align*}\tag{\PointingHand}\T^n1_{[\Psi_n\in x+U]}\le
   \tfrac{\G}{n^{\frac{d}2}}\cdot(e^{-\frac{\eta\|x\|
   ^2}n}+\tfrac1n)\ \ \ 
   \forall\ n\ge 1,\ x\in\Bbb R^d.
 \end{align*}\endproclaim
 \demo{Proof}\ \ We'll show that 
   \begin{align*}\tag{{\scriptsize\faThumbsOUp}}\T^n1_{[\Psi_n\in x+U]}\le
   \tfrac{\G}{n^{\frac{d}2}}\cdot(\frak b(\tfrac{\|x\|}{\sqrt n})e^{-\frac{\eta\|x\|
   ^2}n}+\tfrac1n)\ \ \ 
   \forall\ n\ge 1,\ x\in\Bbb R^d
 \end{align*}
  where 
 $\frak b:\Bbb R\to\Bbb R$ is  a polynomial of degree $3$ with non-negative coefficients and $\frak b(0)=1$.
 \
 
This suffices for (\PointingHand) with 
 suitably adjusted constants $\Gamma>0,\ \eta>0$.
 \

 As before, consider the operators $P(x)\in\text{\tt Hom}\,(L,L)\ \ (x\in\widehat{\Bbb H})$ defined by
 
  \begin{align*}P(x)f=\T(\exp[i\<x,\Psi\>]f).
  \end{align*}

As shown in \cite{HH}, since $E(\|\Psi\|^4)<\infty$,  the function 
$P:\widehat{\Bbb H}\to\text{\tt Hom}\,(L,L)$ is $4$-times continuously differentiable with
 \begin{align*}((\tfrac{\bdy^{k_1}}{\bdy x_{j_1}}\tfrac{\bdy^{k_2}}{\bdy x_{j_2}}&\dots\tfrac{\bdy^{k_q}}{\bdy x_{j_q}})P(x))(f)=i^p\T(\prod_{i=1}^q\Psi_{j_i}^{k_i}\exp[i\<x,\Psi\>]f)
\\ &\ \forall\ 1\le j_1,j_2,\dots,j_q\le d,\ k_1,k_2,\dots,k_q,\ \sum_i k_i=p,\ 1\le p\le 4.  \end{align*}

Moreover, by Nagaev's theorem,\ there are constants $\e>0,\ K>0$ and $\th\in (0,1)$; and  continuous functions
$\l: B(0,\e)\to B_{\Bbb C}(0,1),\  N:B(0,\e)\to \text{\tt Hom}(L,L)$ such that for $x\in  B(0,\e)$,
 \begin{align*}\tag{\dscommercial}\label{dscommercial}  \|P(x)^nf-\l(x)^n N(x)f\|_L\le K\th^n\|f\|_L
 \ \ \ \ \forall\ \ n\ge 1,\ f\in L  
 \end{align*}

where $N(x)$ is a projection onto a one-dimensional subspace (spanned by $N(x)\mathbb{1}$).

\

It follows that $x\mapsto N(x)$ is also $C^4$
 $B(0,\e)\to \text{\tt Hom}(L,L)$ as is $x\mapsto \l(x)\ \ (B(0,\e)\to B_{\Bbb C}(0,1))$.
 \
 
 In particular,  there is a symmetric, positive definite matrix $\g$, and $M>0$ so that for $x\in B(0,\e)$,
 \begin{align*}\tag{\sun}\label{sun}
  \text{\tt Log}\,\l(x)=-\<\g x,x\>(1+\mathcal E(x))\ \text{where}\  |\mathcal E(x)|\le M\|x\|.
 \end{align*}

 Next, fix $h\in L^1(\Bbb H)_+$ so that $1_U\le h$ and $\widehat h\in C^\infty_C(\widehat{\Bbb H})$.
 \footnote{e.g. $h=\widehat g$ with $g(x)=\int_\Bbb H f(y)f(y-x)dm_{\widehat{\Bbb H}}(y)$ for some $f\in C^\infty_C(\widehat{\Bbb H})$ with $\widehat f\ge 1_U$.}
 \
 
 Note that if $h_z(y)=h(y-z)$, then $1_{U+z}\le h_z$ and $\widehat{h_z}(x)=e^{i\<x,z\>}\widehat{h}(x)$ has the same bounded support as $\widehat{h}(x)$.
 \

Possibly increasing $\th\in (0,1)\ \&\ K>0$ in  (\dscommercial) (p. \pageref{dscommercial}), we see by aperiodicity, that 
   \begin{align*} \tag{\dsheraldical}\label{dsheraldical} \sup_{x\in \text{\tt\tiny supp}\,\widehat{h}\setminus B(0,\d)}\,\|P(x)^n\mathbb{1}\|_L\le K\vartheta^n.
   \end{align*}

 \
 
 Thus
   \begin{align*} n^{\frac{d}2}\T^n1_{[\Psi_n\in z+U]}&\le
  n^{\frac{d}2}\T^n(h_z(\Psi_n))\\ & =
  n^{\frac{d}2}\int_{\widehat{\Bbb H}}\widehat{h_z}(x)\T^n(\exp[i\<x,\Psi_n\>])dx\\ & =
  n^{\frac{d}2}\int_{\widehat{\Bbb H}}\widehat{h_z}(x)P(x)^n\mathbb{1}dx\\ & =
  n^{\frac{d}2}\int_{B(0,\d)}\widehat{h_z}(x)\l(x)^ng(x)dx+\mathcal E_n   \end{align*}
  where $g(x):=N(x)\mathbb{1}$ and 
  \begin{align*}  \tfrac{\mathcal E_n}{n^{\frac{d}2}} &=\int_{B(0,\d)}\widehat{h_z}(x)(P(x)^n\mathbb{1}-\l(x)^ng(x))dx
 + \int_{\text{\tt\tiny supp}\,\widehat{h}\setminus B(0,\d)}\widehat{h_z}(x)P(x)^n\mathbb{1}dx,
  \end{align*}
whence by (\dscommercial) and (\dsheraldical) (p. \pageref{dsheraldical})
$$|\mathcal E_n|\le K(1+m_\Bbb H(\text{\tt supp}\,\widehat{h})\|\widehat{h_z}\|_\infty)n^{\frac{d}2}\th^n\le K(1+m_\Bbb H(\text{\tt supp}\,\widehat{h})\|h\|_1)n^{\frac{d}2}\th^n. $$  
  \
  
  Next,  changing variables,
   \begin{align*}\tag{\Bleech}\label{Bleech}
   n^{\frac{d}2}&\int_{B(0,\d)}\widehat{h_z}(x)\l(x)^ng(x)dx\\ &=
   \int_{B(0,\sqrt n \d)}\widehat{h_z}(\tfrac{x}{\sqrt n})\l(\tfrac{x}{\sqrt n})^ng(\tfrac{x}{\sqrt n})dx.\end{align*}
   
     \f{\bf Estimation of {\tt RHS} of (\Bleech)}
    \
    
    Since $\g$ in  (\sun) (p. \pageref{sun}) is positive definite, $\exists\ \eta>0$ so that
    \begin{align*}\tag{{\scriptsize\faTaxi}}\label{faTaxi} \<\g x,x\>\ge\eta\|x\|^2\ \forall\ x\in\Bbb R^d,\end{align*} whence  
$$|\l(\tfrac{x}{\sqrt n})|,\ \exp[-\tfrac{\<\g x,x\>}n]\le e^{-\tfrac{2\eta\|x\|^2}n}\ \forall\ x\in B(0,\sqrt n \d).$$       
       For $n\ge 4,\ x\in B(0,\sqrt n \d)$, set 
     $$\l=\l(\tfrac{x}{\sqrt n})\ \&\ \phi=\exp[-\tfrac{\<\g x,x\>}n],$$
     then $\phi\in (0,1),\ \l\in B(0,1)$ and some baby algebra shows that

 \begin{align*}\tag{\dsmedical}\label{dsmedical} 
 \l^n=\phi^n+
        n\phi^{n-1}(\l-\phi)+(\l-\phi)^2B_n        
        \end{align*}
        where $B_n=B_n(x):=\sum_{k=1}^{n-1}\sum_{j=0}^{k-1}\phi^{n-1-k}\l^j\phi^{k-1-j}$, whence
     \begin{align*}\tag{\dsagricultural}\label{dsagricultural}
      |B_n|&\le \sum_{k=1}^{n-1}|\phi^{n-1-k}\sum_{j=0}^{k-1}\l^j\phi^{k-1-j}|\\ &=
      \sum_{k=1}^{n-1}\sum_{j=0}^{k-1}\phi^{n-2+j}|\l^j|\\ &\le
      n^2e^{-\tfrac{2(n-2)\eta\|x\|^2}n}\le n^2e^{-\eta\|x\|^2}
     \end{align*} 
We now estimate the {\tt RHS} of (\Bleech) according to (\dsmedical) (as on p. \pageref{dsmedical}), 
     \begin{align*}
   \int_{B(0,\sqrt n \d)}&\widehat{h_z}(\tfrac{x}{\sqrt n})\l(\tfrac{x}{\sqrt n})^ng(\tfrac{x}{\sqrt n})dx =\mathcal M_n
  +\mathcal Q_n+\mathcal D_n\\ &=:\int_{B(0,\sqrt n \d)}\widehat{h_z}(\tfrac{x}{\sqrt n})g(\tfrac{x}{\sqrt n})\cdot\(\phi^n+
        n\phi^{n-1}(\l-\phi)+(\l-\phi)^2B_n\)\cdot dx
  \end{align*}
   where
      $$\mathcal M_n:=\int_{B(0,\sqrt n \d)}e^{i\<x,\tfrac{z}{\sqrt n}\>}
      \widehat{h}(\tfrac{x}{\sqrt n})g(\tfrac{x}{\sqrt n})\exp[-\<\g x,x\>]dx$$
      is the ``{\tt main term}'' with ``{\tt error terms}''

   $$\mathcal Q_n=n\int_{B(0,\sqrt n \d)}e^{i\<x,\tfrac{z}{\sqrt n}\>}\widehat{h}(\tfrac{x}{\sqrt n})g(\tfrac{x}{\sqrt n})\exp[-\<\g x,x\>\cdot\tfrac{n-1}n](\l(\tfrac{x}{\sqrt n})-\exp[-\tfrac{\<\g x,x\>}n])dx$$
   and
   $$\mathcal D_n=\int_{B(0,\sqrt n \d)}e^{i\<x,\tfrac{z}{\sqrt n}\>}\widehat{h}(\tfrac{x}{\sqrt n})g(\tfrac{x}{\sqrt n})\exp[-\<\g x,x\>\cdot\tfrac{n-2}n](\l(\tfrac{x}{\sqrt n})-\exp[-\tfrac{\<\g x,x\>}n])^2B_n(x)dx.$$

There is a polynomial $\frak q(x)=a_0+a_1(x)+a_2(x)+a_3(x)$ (where each $a_i$ is a sum of monomials of degree $i$) so that
$$\widehat{h}(y)g(y)=\frak q(y)+O(\|y\|^4),$$ whence
 \begin{align*}\tag{\HandWash}\label{HandWash}
   \widehat{h}(\tfrac{x}{\sqrt n})g(\tfrac{x}{\sqrt n})=\frak q(\tfrac{x}{\sqrt n})+O(\tfrac{\|x\|^4}{n^2}).
 \end{align*}

\f{\bf The  term $\mathcal M_n$} 
\

The following estimations are uniform in $z\in\Bbb R^d$.
\

By (\HandWash), we have 
$$\mathcal M_n =\int_{B(0,\sqrt n \d)}e^{i\<x,\tfrac{z}{\sqrt n}\>}\frak q(\tfrac{x}{\sqrt n})\exp[-\<\g x,x\>]dx +O(\tfrac1{n^2}).$$
For some $\vartheta\in (0,1)$ we have
 \begin{align*}
  \tag{\Dontwash}\label{Dontwash}
\int_{\Bbb R^d\setminus B(0,\sqrt n \d)}\|x\|^j\exp[-\<\g x,x\>]dx=O(\vartheta^n)\ \ (0\le j\le 4).
 \end{align*}

Thus
 \begin{align*}\mathcal M_n=
 \int_{\Bbb R^d}e^{i\<x,\tfrac{z}{\sqrt n}\>}\frak q(\tfrac{x}{\sqrt n})\exp[-\<\g x,x\>]dx+O(\tfrac1{n^2})
   \end{align*}

and
 \begin{align*}
  \int_{\Bbb R^d}e^{i\<x,\tfrac{z}{\sqrt n}\>}\frak q(\tfrac{x}{\sqrt n})\exp[-\<\g x,x\>]dx&=\sum_{j=0}^3
  \tfrac1{n^{\frac{j}2}}\int_{\Bbb R^d}e^{i\<x,\tfrac{z}{\sqrt n}\>}a_j(x)\exp[-\<\g x,x\>]dx\\ &=
  \sum_{j=0}^3\tfrac1{n^{\frac{j}2}}(a_j(\nabla)\phi)(\tfrac{z}{\sqrt n})\\ &=
  \sum_{j=0}^3\tfrac1{n^{\frac{j}2}}b_j(\tfrac{z}{\sqrt n})\phi(\tfrac{z}{\sqrt n})
   \end{align*} 
where $b_j$ is a sum of degree $j$ monomials,  $b_0=1$ and $\phi(y)=\int_{\Bbb R^d}e^{i\<x,y\>}\exp[-\<\g x,x\>]dx$.

Evidently, $\phi(y)=c e^{-\<\g y,y\>}\le ce^{-\eta\|y\|^2} $ where $c>0$ and $\eta>0$ is as in 
({\scriptsize\faTaxi}) (as on p. \pageref{faTaxi}).
\

Moreover, 
$$|\sum_{j=0}^3\tfrac1{n^{\frac{j}2}}b_j(\tfrac{x}{\sqrt n})|\le \frak b(\tfrac{\|x\|}{\sqrt n})\ \forall\ x\in\Bbb R^d$$
where
 $\frak b:\Bbb R\to\Bbb R$ is  a polynomial of degree $3$ with non-negative coefficients and $\frak b(0)=1$.

It follows that, uniformly in $z$:
$$\mathcal M_n\le \frak b(\tfrac{z}{\sqrt n})e^{-\frac{\eta\|z\|^2}n}+O(\tfrac1{n^2}).$$
                                                                  
   \f{\bf The  term $\mathcal D_n$}
        \
        
   To estimate $\mathcal D_n$ we use (\dsagricultural) as on p. \pageref{dsagricultural} and
    \begin{align*}\tag{\dsarchitectural}\label{dsarchitectural}
        |\l(\tfrac{x}{\sqrt n})-\exp[-\tfrac{\<\g x,x\>}n]|\le A\tfrac{\|x\|^3}{n^{\frac32}}
    \end{align*}
obtaining
    \begin{align*}   |\mathcal D_n|&\le 
    \int_{B(0,\sqrt n \d)}\exp[-\<\g x,x\>\cdot\tfrac{n-2}n]|\l(\tfrac{x}{\sqrt n})-\exp[-\tfrac{\<\g x,x\>}n]|^2|B_n(x)|dx\\
    & \le \frac{A^2}{n^3}\int_{B(0,\sqrt n \d)}\|x\|^6|B_n(x)|\exp[-\tfrac{\<\g x,x\>}2]dx\ \text{by (\dsarchitectural)} \\
    & \le \frac{A^2}{n}\int_{B(0,\sqrt n \d)}\|x\|^6\exp[-\tfrac{\<\g x,x\>}2]dx\ \text{by (\dsagricultural)}.    
    \end{align*}
       \f{\bf The  term $\mathcal Q_n$}
       \
       
This needs a  Taylor expansion of higher order than (\dsarchitectural) (p. \pageref{dsarchitectural}):
\

For $n\ge 1$ and $x\in B(0,\sqrt n \d)$,
    \begin{align*}
    \l(\tfrac{x}{\sqrt n})-\exp[-\tfrac{\<\g x,x\>}n]=\tfrac1{n^{\frac32}}c_3(x)+\frak e_n(x)     
    \end{align*}
where $c_3$ is a sum of degree $3$ monomials and $|\frak e_n(x)|\le \frac{C\|x\|^4}{n^2}$.
\

Thus 
\begin{align*} 
     \mathcal Q_n&=
     n\int_{B(0,\sqrt n \d)}e^{i\<x,\tfrac{z}{\sqrt n}\>}
     \widehat{h}(\tfrac{x}{\sqrt n})g(\tfrac{x}{\sqrt n})\exp[-\<\g x,x\>\cdot\tfrac{n-1}n]
     (\l(\tfrac{x}{\sqrt n})-\exp[-\tfrac{\<\g x,x\>}n])dx\\ &=
     n\left[\int_{B(0,\sqrt n \d)}e^{i\<x,\tfrac{z}{\sqrt n}\>}\frak q(\tfrac{x}{\sqrt n})\exp[-\<\g x,x\>\cdot\tfrac{n-1}n]
     (\l(\tfrac{x}{\sqrt n})-\exp[-\tfrac{\<\g x,x\>}n])dx+O(\frac1{n^2})\right]\\ &
     \ \ \ \ \ \ \ \ \ \ \ \ \ \ \text{by (\HandWash) on p. \pageref{HandWash}}\\ &=n\int_{B(0,\sqrt n \d)}e^{i\<x,\tfrac{z}{\sqrt n}\>}\frak q(\tfrac{x}{\sqrt n})\exp[-\<\g x,x\>\cdot\tfrac{n-1}n]
     (\l(\tfrac{x}{\sqrt n})-\exp[-\tfrac{\<\g x,x\>}n])dx+O(\frac1{n})\\ &=
   n\int_{B(0,\sqrt n \d)}e^{i\<x,\tfrac{z}{\sqrt n}\>}\frak q(\tfrac{x}{\sqrt n})\exp[-\<\g x,x\>\cdot\tfrac{n-1}n]
     (\tfrac1{n^{\frac32}}c_3(x)+\frak e_n(x) )dx+O(\frac1{n})\\ &=
     \tfrac1{\sqrt n}\int_{B(0,\sqrt n \d)}e^{i\<x,\tfrac{z}{\sqrt n}\>}\frak q(\tfrac{x}{\sqrt n})\exp[-\<\g x,x\>\cdot\tfrac{n-1}n]
     c_3(x)dx\\ &\ \ \ \ \ \ \ \ \ \ \ \ \ \ \ \    +\int_{B(0,\sqrt n \d)}e^{i\<x,\tfrac{z}{\sqrt n}\>}\frak q(\tfrac{x}{\sqrt n})\exp[-\<\g x,x\>\cdot\tfrac{n-1}n]\frak e_n(x)dx+O(\frac1{n})\\ &=
     \tfrac1{\sqrt n}\int_{B(0,\sqrt n \d)}e^{i\<x,\tfrac{z}{\sqrt n}\>}\frak q(\tfrac{x}{\sqrt n})\exp[-\<\g x,x\>\cdot\tfrac{n-1}n]
     c_3(x)dx+O(\frac1{n}).
\end{align*}
Now 
   \begin{align*}\int_{B(0,\sqrt n \d)}&e^{i\<x,\tfrac{z}{\sqrt n}\>}\frak q(\tfrac{x}{\sqrt n})\exp[-\<\g x,x\>\cdot\tfrac{n-1}n]
     c_3(x)dx\\ &=\int_{\Bbb R^d}e^{i\<x,\tfrac{z}{\sqrt n}\>}\frak q(\tfrac{x}{\sqrt n})\exp[-\<\g x,x\>\cdot\tfrac{n-1}n]
     c_3(x)dx+O(\vartheta^n)\ \ \text{by (\Dontwash) on p. \pageref{Dontwash}}\\ &=
     \sum_{j=0}^3\frac1{n^{\frac{j}2}}\int_{\Bbb R^d}a_j(x)c_3(x)e^{i\<x,\tfrac{z}{\sqrt n}\>}\frak q(\tfrac{x}{\sqrt n})\exp[-\<\g x,x\>\cdot\tfrac{n-1}n]
     dx+O(\vartheta^n)\\ &=\sum_{j=0}^3\frac1{n^{\frac{j}2}}(a_j(\D)c_3(\D)\phi^{(n)})(\tfrac{z}{\sqrt n})\ \ \text{where}\ 
     \phi^{(n)}(u)=\int_{\Bbb R^d}\exp[i\<u,x\>-\<\g x,x\>\cdot\tfrac{n-1}n]dx,
    \\ &=\sum_{j=0}^3\frac1{n^{\frac{j}2}}\frak p_j(\tfrac{z}{\sqrt n})\phi(\tfrac{z}{\sqrt n})
       \end{align*} where the  $\frak p_j$ are polynomials.\ \ \ \CheckedBox\ ($\diamondsuit$)
       
       \
       
       This proves ({\scriptsize\faThumbsOUp}) and hence (\PointingHand).\ \Checkedbox
      
      \
       
       \
       
     \proclaim{Deviation Lemma}
 \

 {\rm (i)}\ There are constants $\G,\ \eta>0$ 
 so that $\forall\ M>0\ \exists\ t_0(M)>0$ such that
  \begin{align*}\T^n1_{A\cap [\frak r_n\in I+t-y]}&\le 
   \G\(\frac{\exp[-\frac{\eta |n-\l t|^2}t]}{\sqrt t}
   +\frac1{n^{\frac32}}\)\\ &
   \forall\ M>0,\ t>t_0(M)\ \&\ n\in B_K(M,t)\end{align*}
    where $K>1\ \&\ B_K(M,t)$ are as on p. \pageref{dsrailways}.
    \
    
{\rm (ii)}\ If in addition, $E(\|\v\|^4)<\infty$, then,  possibly changing $\G,\ \eta\ \&\ t_0(M)$, 
we have for $M>0,\ t>t_0(M),\ x(t)\in\Bbb R^\kappa,\ 
\|x(t)\|\le \tfrac{M\sqrt t}{2\l K^2}$:

\begin{align*}\T^n&1_{A\cap [\varphi_n\in  x(t)-z+V,\ \frak r_n\in I+t-y]}\le
   \frac{\G}{t^{\frac{\kappa}2}}\cdot\(\frac{\exp[-\frac{\eta|n-\l t|^2}t]}{\sqrt t}+\frac1{n^{\frac32}}\)\\ &
   \forall\  n\in B_K(M,t).\end{align*}\endproclaim
   \demo{Proof of (ii) in the deviation lemma}
   \
   
   Let $\psi:=(\v,\frak r-\varkappa):X\to\Bbb G\x\Bbb R$ where $\varkappa:=E(\frak r)=\tfrac1\l$. 
   
   For $t>0,\ x(t)\in\Bbb R^\kappa,\ \|x(t)\|\le \tfrac{M\sqrt t}{4\l K^2},\ n\in  B_K(M,t)$,
   \begin{align*}\T^n1_{A\cap [\varphi_n\in  x(t)-z+V,\ \frak r_n\in I+t-y]}&=
   \T^n1_{A\cap [\psi_n\in  x_{n,t}+V\x I]}\end{align*}
   where $x_{n,t}:=(x(t)-z,t-n\varkappa-y)$.
   \
   
   By the assumptions, the extended {\tt LLT} estimate applies and it follows from (\PointingHand) that
   \ there exist
 $\Gamma>0,\ \eta>0$ 
 so that
   \begin{align*}\T^n1_{[\psi_n\in x+U]}\le
   \frac{\G}{n^{\frac{\kappa}2}}\cdot(\frac{\exp[-\frac{\eta\|x_{n,t}\|^2}n]}{\sqrt n}+\frac{\G}{n^{\frac32}})\ \ \ 
   \forall\ n\in B_K(M,t).\end{align*}
             Next \begin{align*}\|x_{n,t}\|&\ge |t-n\varkappa|-\|x(t)\|-(|y|+\|z\|)\\ &=
             \tfrac{|n-\l t|}\l-\|x(t)\|-(|y|+\|z\|)\\ &=\tfrac{|n-\l t|}{2\l}+\tfrac{|n-\l t|}{2\l}-\|x(t)\|-(|y|+\|z\|).   \end{align*} 
             Now
             \begin{align*}\tfrac{|n-\l t|}{2\l}&-\|x(t)\|-(|y|+\|z\|)\ge
             \tfrac{M\sqrt n}{2\l K}-\tfrac{M\sqrt t}{2\l K^2}-(|y|+\|z\|)\\ &\ge
             \tfrac{M\sqrt t}{2\l K^2}-(|y|+\|z\|)>0 \ \forall\ M>0,\ t>t_0(M):=\tfrac{4\l^2K^2(|y|+\|z\|)^2}{M^2}.
                           \end{align*}

       Setting $y_{n,t}:=\tfrac{x_{n,t}}{\sqrt n}$, we have for $M>0,\ t>t_0(M)$ and $n\in B_K(M,t)$,
       \begin{align*}\|y_{n,t}\|>\tfrac{|n-\l t|}{2\l\sqrt n}.
       \end{align*}

As before, let $\xi>0$ be so that
    $$\<\g x,x\>\ge \xi\|x\|^2\ \forall\ x\in\Bbb R^{\kappa+1}.$$
It follows that  for $M>0,\ t>t_0(M)$ and $n\in B_K(M,t)$,
  $$\tfrac{\<\g x_{n,t},x_{n,t}\>}n\ge \<\g y_{n,t},y_{n,t}\>\ge \xi\|y_{n,t}\|^2\ge \eta\tfrac{|\l t-n|^2}t$$
  where $\eta:=\tfrac{\xi}{4K\l^2}$.

Thus,
  \begin{align*}\T^n1_{[\psi_n\in x+U]}&\le
   \frac{K^2\G}{t^{\frac{\kappa}2}}\cdot(\tfrac{\exp[-\frac{\eta|n-\l t|^2}t]}{\sqrt t}+\tfrac{\G}{n^{\frac32}})\\ &
   \forall\ M>0,\ t>t_0(M)\ \&\ n\in B_K(M,t).\ \CheckedBox\end{align*}
   
\end{document}